\documentclass[a4paper,11pt]{amsart}
\usepackage{amsmath}
\usepackage{amssymb}
\usepackage{amsthm}
\usepackage[all]{xy}
\usepackage[latin1]{inputenc}        % accents
\usepackage[dvips]{graphics}
\usepackage[dvips]{graphicx}
\usepackage{amscd}
\usepackage{mathrsfs}
\usepackage[mathcal]{eucal}
\usepackage{fullpage}
\author{Francesco Polizzi}
\address{Dipartimento di Matematica, Università della Calabria, Via Pietro Bucci,
87036 Arcavacata di Rende (CS), Italy.}
\email{polizzi@mat.unical.it}
\title[Surfaces isogenous to a product]
{On surfaces of general type with $p_g=q=1$ isogenous to a product
of curves}
\date{\today}

\newtheorem{inizio}{Lemma}[section]
\newtheorem{theorem}[inizio]{Theorem}
\newtheorem{corollary}[inizio]{Corollary}
\newtheorem{proposition}[inizio]{Proposition}
\newtheorem{lemma}[inizio]{Lemma}
\newtheorem{definition}[inizio]{Definition}

\newtheorem*{teoA}{Theorem A}
\newtheorem*{teoB}{Theorem B}
\newtheorem*{teoC}{Theorem C}
\newtheorem*{prob}{Main Problem}
\theoremstyle{definition}
\newtheorem{remark}[inizio]{Remark}

\newcommand{\lr}{\longrightarrow}
\newcommand{\hG}{\Gamma}

\newcommand{\mO}{\mathcal{O}}

\begin{document}
\subjclass{14J29 (primary), 14J10, 20F65} \keywords{Surfaces of
general type, actions of finite groups on curves}
\maketitle
\begin{center}
\emph{To the memory of my colleague and friend Giulio Minervini.}
\end{center}
\abstract A smooth algebraic surface $S$ is said to be  \emph{isogenous to a product of unmixed type} if there exist two smooth curves $C,
\;F$ and a finite group $G$, acting faithfully on both $C$ and $F$
and freely on their product, so that $S=(C \times F)/G$. In this
paper we classify the surfaces of general type with $p_g=q=1$ which
are isogenous to an unmixed product, assuming that the group $G$ is
abelian. It turns out that they belong to four families, that we
call surfaces of type $I, \;II, \; III, \; IV$. The moduli spaces
$\mathfrak{M}_{I}, \; \mathfrak{M}_{II}, \; \mathfrak{M}_{IV}$ are
irreducible, whereas $\mathfrak{M}_{III}$ is the disjoint union of
two irreducible components. In the last section
 we start the analysis of the case where $G$
is not abelian, by constructing several examples.
\endabstract
% \tableofcontents
\section{Introduction}
The problem of classification of surfaces of general type is of
exponential computational complexity, see \cite{Ca92}, \cite{Ch96},
\cite{Man97}; nevertheless, one can hope to classify at least those
with small numerical invariants. It is well-known that the first
example of surface of general type with $p_g=q=0$ was given by
Godeaux in \cite{Go31}; later on, many other examples were
discovered. On the other hand, any surface $S$ of general type
verifies $\chi(\mathcal{O}_S)>0$, hence $q(S)>0$ implies $p_g(S)>0$.
It follows that the surfaces with $p_g=q=1$ are the irregular ones
with the lowest geometric genus, hence it would be important to
achieve their complete classification; so far, this has been
obtained only in the cases $K_S^2=2,3$ (see \cite{Ca}, \cite{CaCi1},
\cite{CaCi2}, \cite{Pol05}, \cite{CaPi05}). As the title suggests,
this paper considers surfaces of general type with $p_g=q=1$ which
are \emph{isogenous to a product}. This means that there exist two
smooth curves $C, \; F$ and a finite group $G$, acting freely on
their product, so that $S=(C \times F)/G$. We have two cases: the
\emph{mixed} case, where the action of $G$ exchanges the two factors
(and then $C$ and $F$ are isomorphic) and the \emph{unmixed} case,
where $G$ acts diagonally. In the unmixed case  $G$ acts
separately on $C$ and $F$, and the two projections $\pi_C \colon C
\times F \lr C, \; \pi_F \colon C \times F \lr F$ induce two
isotrivial fibrations $\alpha \colon S \lr C/G, \; \beta \colon S
\lr F/G$, whose smooth fibres are isomorphic to $F$ and $C$,
respectively.  If $S$ is isogenous to a product, there exists a
unique realization $S=(C \times F) /G$ such that the genera $g(C),
\; g(F)$ are minimal (\cite{Ca00}, Proposition 3.13); we will always work
with minimal realizations. Surfaces of general type with $p_g=q=0$
 isogenous to a product appear in \cite{Be2}, \cite{Par03} and \cite{BaCa03}; their
 complete classification has been finally obtained in \cite{BaCaGr06}.
Some unmixed examples with $p_g=q=1$ have been given in
\cite{Pol06}; so it seemed natural to attack the following
\begin{prob}
Classify all surfaces of general type with $p_g=q=1$ isogenous to a
product, and describe the corresponding irreducible components of
the moduli space.
\end{prob}
In this paper we fully solve the Main Problem in the unmixed case
assuming that
 the group $G$ is abelian.
Our results are the following:
\begin{teoA} [see Theorem \ref{isog prod}]
If  the group $G$ is abelian, then there exist exactly four
families of surfaces of general type with $p_g=q=1$ isogenous to an
unmixed product. In every case $g(F)=3$, whereas the occurrences for
$g(C)$ and $G$ are
\begin{itemize}
\item[$I.$] $g(C)=3, \quad G=(\mathbb{Z}_2)^2$;
\item[$II.$] $g(C)=5, \quad G=(\mathbb{Z}_2)^3$;
\item[$III.$] $g(C)=5, \quad G=\mathbb{Z}_2 \times \mathbb{Z}_4$;
\item[$IV.$] $g(C)=9, \quad G= \mathbb{Z}_2 \times \mathbb{Z}_8$.
\end{itemize}
Surfaces of type $I$ already appear in  \cite{Pol06}, whereas those
of type $II, \; III, \; IV$ provide new examples of minimal surfaces
of general type with $p_g=q=1, \; K^2=8$.
\end{teoA}

\begin{teoB} [see Theorem \ref{modulispaces}]
The moduli spaces $\mathfrak{M}_I, \; \mathfrak{M}_{II}, \;
\mathfrak{M}_{IV}$ are irreducible of dimension $5,\;4,\;2$,
respectively. The moduli space $\mathfrak{M}_{III}$ is the disjoint
union of two irreducible components $\mathfrak{M}_{III}^{(1)},
\;\mathfrak{M}_{III}^{(2)}$, both of dimension $3$.
\end{teoB}

The case where $G$ is not abelian
is more difficult, and a complete classification is still lacking
(see Remark \ref{achieve-non-ab}). However, we can shed some light
on this problem, by proving

\begin{teoC}[see Theorem \ref{noabeliano}]
Let $S=(C \times F) /G$ be a surface of general type with $p_g=q=1$,
isogenous to an unmixed product, and assume that the group $G$ is
not abelian. Then the following cases occur:
\begin{equation*}
\begin{tabular}{c|c|c|c}
$G$ & $|G|$ & $g(C)$ & $g(F)$ \\
 \hline
 $S_3$ &  $6$ & $3$ & $4$  \\
 $D_4$ & $8$ & $3$ & $5$ \\
 $D_6$ & $12$ & $7$ & $3$  \\
 $A_4$ & $12$ & $4$ & $5$  \\
 $S_4$ & $24$ & $9$ & $4$ \\
 $A_5$ & $60$ & $21$ & $4$ 
\end{tabular}
\end{equation*}
\end{teoC}
The examples with $G= S_3$ and $D_4$ already appear
 in \cite{Pol06}, whereas the others are new. It would be
interesting to have a description of the moduli spaces for these
new examples (see Remark \ref{moduli-new}). \\
While describing the organization of the paper we shall now explain
the steps of our classification procedure in more detail. The crucial point is that in the unmixed case the geometry of the surface $S=(C \times F)/G$
is encoded in the geometry of the two $G-$covers $h \colon C \lr C/G$, $f
\colon F \lr F/G$. This allows us to "detopologize" the problem by 
transforming it into an equivalent problem about the existence of a pair
of epimorphisms from two groups of Fuchsian type into $G$; this
is essentially an application of the Riemann's existence theorem. 
These epimorphisms must satisfy some additional properties in order to get a free action of $G$ on $C \times F$ and a quotient surface with the desired invariants (Proposition \ref{structureresult1}). The
geometry of the moduli spaces can be also recovered from these algebraic
data (Propositions \ref{structuremoduli} and \ref{structuremoduli1}).
\\
In the nonabelian case we follow a similar approach (Proposition \ref{structureresult2}). \\
In Section \ref{topological} we fix the algebraic set up. The reader that
is only interested in the proof of Theorems $A$ and $C$
might skip to Section \ref{surfisog} after reading Section \ref{admiss-epi}.
On the other hand, the content of Sections \ref{hurwitz-moves}, \ref{g'=0}, \ref{g'=1,n=1}, \ref{g'=1,n=2}  is essential in order to understand the proof of Theorem $B$. The results in \ref{g'=0}
 are well-known, whereas for those in \ref{g'=1,n=1} and \ref{g'=1,n=2}
 we have not been able to find any complete reference; so we had
 to carry out "by hand" all the (easy) computations. \\
 In Section \ref{surfisog} we establish some basic results about surfaces
 $S$ of general type with $p_g=q=1$ isogenous to a product. Such surfaces
 are always minimal and verify $K_S^2=8$. Moreover, we show that if $G$ is
 abelian then the Albanese fibration of $S$ is a genus $3$ pencil with two
 double fibres. \\ 
The main results of Section \ref{building} are Propositions \ref{structureresult1} and \ref{structuremoduli1},
which play a central role in this paper as they translate our Main Problem
"from geometry to algebra". \\  
Section \ref{classification} contains the proof of Theorem $A$,
whereas Section \ref{secmodulispaces} contains the proof of Theorem $B$.  \\
In Section \ref{paracansystem}  we study the paracanonical system
$\{K\}$ for surfaces of type $I, \; II, \; III, \; IV$, showing that
in any case it has index $1$ (Theorem \ref{calcolo indice}). This
section could appear as a digression with respect to the main
 theme of the paper; however, since the index of $\{K\}$ is an important
 invariant of $S$ (see \cite{Be88}, \cite{CaCi1}, \cite{CaCi2}) we thought worthwhile computing it. \\
Finally, Section \ref{newexamples} deals with the proof of Theorem $C$.\\ \\
$\mathbf{Notations\;\;and\;\;conventions.}$ All varieties and
morphisms in this article are defined over the field $\mathbb{C}$ of
complex numbers. By ``surface'' we mean a projective, non-singular
surface $S$, and for such a surface $K_S$ or $\omega_S$ denote the
canonical class, $p_g(S)=h^0(S, \; K_S)$ is the \emph{geometric
genus}, $q(S)=h^1(S, \; K_S)$ is the \emph{irregularity} and
$\chi(\mathcal{O}_S)=1-q(S)+p_g(S)$ is the \emph{Euler
characteristic}. If $S$ is a surface with $p_g=q=1$, then $\alpha
\colon S \longrightarrow E$ is the Albanese map of $S$ and $F$
denotes the general fibre of $\alpha$.   \\ \\
$\mathbf{Acknowledgments.}$ Part of the present research was done
during the author's visits to the Universities of Bayreuth,
Cambridge and Warwick, where he was supported by Eu Research
Training Network EAGER, no. HPRN-CT-2000-00099, and by a Marie Curie
pre-doc fellowship. The author is grateful to C. Ciliberto for
useful advice and constructive remarks during the preparation of
this work and to  F. Catanese and I. Bauer for
 suggesting the problem and sharing the ideas contained in
\cite{BaCa03}. He also wishes to thank  A. Corti, M. Reid, G. Infante for helpful conversations and R. Pardini for pointing out
some mistakes contained in the first version of this paper. Finally, he is
indebted with the referee for several detailed comments that
considerably improved the presentation of these results.

\section{Topological background} \label{topological}
Many of the result that we collect in this section are standard, so
proofs are often omitted. We refer the reader to [Br90, Section 2],
[Bre00, Chapter 3] and \cite{H71} for more details.
\subsection{Admissible epimorphisms} \label{admiss-epi}
Let us denote by $\Gamma=\Gamma(g' \;|\; m_1, \ldots, m_r)$ the
abstract group of Fuchsian type with a presentation of the form
\begin{equation}  \label{presfuchs}
\begin{split}
\textrm{generators:} & \,\,\,\,\, a_1, \ldots, a_{g'}, \; b_1, \ldots,b_{g'}, \; c_1, \ldots ,c_r \\
\textrm{relations:} & \,\,\,\,\, c_1^{m_1}= \cdots = c_r^{m_r}=1 \\
 & \,\,\,\,\, c_1c_2 \cdots c_r \Pi_{i=1}^{g'}[a_i,b_i]=1.
\end{split}
\end{equation}
The \emph{signature} of $\hG$ is the ordered set of integers $(g'\;
| \; m_1, \ldots,m_r)$, where without loss of generality we may
suppose $2 \leq m_1 \leq m_2\leq \cdots \leq m_r$. We will call $g'$
the $\emph{orbit genus}$ of $\Gamma$ and $\mathbf{m}:=(m_1, \ldots,
m_r)$ the $\emph{branching data}$. In fact the group  $\Gamma$  acts
on the upper half-plane $\mathscr{H}$ so that the quotient space
$\mathscr{H}/\hG$ is a compact Riemann surface of genus $g'$ and the
$m_i$ are the ramification numbers of the branched covering
$\mathscr{H} \lr \mathscr{H}/\Gamma$. For convenience we make
abbreviations such as $(2^3,3^2)$ for $(2,2,2,3,3)$
 when we write down the branching data.
If the branching data are empty, the corresponding group
$\Gamma(g'\; | \;-)$ is isomorphic to the fundamental group of a compact
 Riemann surface of genus $g'$; it will be denoted by $\Pi_{g'}$. The
following result, which is essentially a reformulation of the Riemann's existence
theorem, translates the problem of finding Riemann
 surfaces with automorphisms into the group theoretic problem of finding
 certain normal subgroups in a given group of Fuchsian type.

\begin{proposition} \label{Riemann-ext}
 A finite group $G$ acts as a group of automorphisms of some compact Riemann
 surface $X$ of genus $g \geq 2$ if and only if there exist a group of Fuchsian
 type $\Gamma=\Gamma(g' \;|\; m_1, \ldots, m_r)$ and an epimorphism $\theta
 \colon \Gamma \lr G$ such that $\emph{Ker}\;\theta\cong \Pi_g$.
\end{proposition}
Since $\Pi_g$ is torsion-free, it follows that $\theta$ preserves
the orders of the elliptic generators $c_1, \ldots , c_r$ of
$\Gamma$. This motives the following

\begin{definition} \label{admepi}
Let $G$ be a finite group. An epimorphism $\theta \colon  \Gamma\lr
G $ is called \emph{admissible} if  $\theta(c_i)$ has order $m_i$
for every $i\in \{ 1, \ldots, r \}$. If an admissible epimorphism $\theta \colon  \Gamma\lr
G$ exists, then $G$ is said to be $(g' \;|\; m_1, \ldots,
m_r)-$\emph{generated}.
\end{definition}
\begin{proposition} \label{ab-no1}
If an abelian group $G$ is $(g'\; |\; m_1, \ldots, m_r)-$generated,
then $r \neq 1$.
\end{proposition}
\begin{proof}
Suppose $G$  abelian and $r=1$. Then relation $x_1 \Pi_{i=1}^{g'}[a_i,b_i]=1$
yields $\theta(x_1)=0$ for any epimorphism $\theta \colon \Gamma \lr G$, so $\theta$ cannot be admissible.
\end{proof}
If $G$ is $(g' \;|\; m_1, \ldots, m_r)-$generated, set
\begin{equation*}
\begin{split}
g_i &:=\theta(c_i) \quad  1 \leq i \leq r; \\
h_j &:=\theta(a_j) \quad  1 \leq j \leq g'; \\
h_{j+g'} &:=\theta(b_j) \quad 1 \leq j \leq g'.
\end{split}
\end{equation*}
The elements $g_1, \ldots, g_r, \; h_1, \ldots ,h_{2g'}$ generate
$G$ and moreover one has
\begin{equation} \label{zz}
 g_1g_2\cdots g_r \Pi_{i=1}^{g'}
[h_i,h_{i+g'}]=1
\end{equation}
and
\begin{equation} \label{zz1}
o(g_i)=m_i.
\end{equation}
\begin{definition}
An \emph{admissible generating vector} $($or, briefly, a
\emph{generating vector}$)$ of $G$ with respect to $\Gamma$ is a
$(2g'+r)$-ple of elements
\begin{equation*}
\mathcal{V}=\{g_1, \ldots, g_r; \; h_1, \ldots \, h_{2g'} \}
\end{equation*}
such that $\mathcal{V}$ generates $G$ and $(\ref{zz})$, $(\ref{zz1})$
are satisfied.
\end{definition}
If $G$ is abelian we use the additive notation and relation
(\ref{zz}) becomes
\begin{equation} \label{zz2}
g_1+ \cdots +g_r=0.
\end{equation}
It is evident that giving a generating vector for $G$ with respect to
$\Gamma$ is equivalent to give an admissible epimorphism $\theta
\colon \Gamma \lr G$; such an epimorphisms fixes the representation
of $G$ as a group of conformal automorphisms of a compact Riemann
surface $X$ of genus $g$ and the quotient $X/G$ has genus $g'$,
where $g$ and $g'$ are related by the Riemann-Hurwitz formula
\begin{equation} \label{riemanhur}
2g-2=|G| \left( 2g'-2+\sum_{i=1}^r\bigg(1-\frac{1}{\;m_i} \bigg)
 \right).
\end{equation}
Hence, accordingly to Proposition \ref{Riemann-ext}, there is a
short exact sequence
\begin{equation} \label{shortth}
1 \lr \Pi_g \stackrel{\iota_{\theta}}{\lr} \Gamma
\stackrel{\theta}{\lr} G \lr 1
\end{equation}
such that $\Gamma$ can be viewed as the \emph{orbifold fundamental
group}
 of the branched cover $X \lr X/G$ (see \cite{Ca00}). In particular, the
 cyclic subgroups $\langle g_i \rangle$ and their conjugates are the non-trivial
 stabilizers of the action of $G$ on $X$. 
 
\subsection{Hurwitz moves} \label{hurwitz-moves} 
 
 Looking at exact sequence (\ref{shortth}) it is important to remark
 that $X$ is defined up to automorphisms not by the specific $\theta$, but
 rather by its kernel $\iota_{\theta}(\Pi_g)$; this motives the following
\begin{definition}
We set
\begin{equation*}
\emph{Epi}(\Pi_g, \Gamma, G):= \left \{ \begin{array}{l}
\emph{Admissible epimorphims\; \;} \theta \colon \Gamma \lr G \\
\emph{such that\;\;} \emph{Ker\;} \theta \cong \Pi_g
\end{array} \right \} / \sim
\end{equation*}
where $\theta_1 \sim \theta_2$ if and only if
$\emph{Ker\;}\theta_1=\emph{Ker\;}\theta_2$.
\end{definition}
Abusing notation we will often not distinguish between an epimorphism $\theta$
 and its class in $\textrm{Epi}(\Pi_g, \Gamma, G)$.
An automorphism $\eta \in \textrm{Aut}(\Gamma)$ is said to be
orientation-preserving if, for all $i \in \{1, \ldots, r\}$, there exists
$j$ such that $\eta(c_i)$ is
conjugated to $c_j$. This of course implies $o(c_i)=o(c_j)$. The
subgroup of orientation-preserving automorphisms of $\Gamma$ is
denoted by $\textrm{Aut}^+(\Gamma)$ and the quotient
$\textrm{Mod}(\Gamma):=\textrm{Aut}^+(\Gamma) /
\textrm{Inn}(\Gamma)$ is called the \emph{mapping class group} of
$\Gamma$. There is a natural action of $\textrm{Aut}(G) \times
\textrm{Mod}(\Gamma)$ on $\textrm{Epi}(\Pi_g, \Gamma, G)$, namely
\begin{equation*}
(\lambda, \eta)\cdot \theta:= \lambda \circ \theta \circ \eta.
\end{equation*}
\begin{proposition} \label{same-class}
Two admissible epimorphisms $\theta_1, \; \theta_2 \in
\emph{Epi}(\Pi_g, \Gamma, G)$ define the same equivalence class of
$G-$actions if and only if they lie in the same  $\emph{Aut}(G)
\times \emph{Mod}(\Gamma)-$class.
\end{proposition}
The non-trivial part of the proof is to show that  $\textrm{Aut}(G)
\times \textrm{Mod}(\Gamma)-$equivalent epimorphisms give equivalent
$G-$actions; this depends on  Teichmüller theory and proofs can be
found in \cite{McB66} and \cite{H71}. \\
The action of $\textrm{Aut}(G) \times \textrm{Mod}(\Gamma)$ on
$\textrm{Epi}(\Pi_g, \Gamma, G)$ naturally induces an action on the
set of generating vectors (up to inner automorphisms of $G$); in
particular, if $\theta_1$ and $\theta_2$ are in the same $\{
\textrm{Id}\} \times \textrm{Mod}(\Gamma)-$class,
 we say that the corresponding generating vectors are
related by a \emph{Hurwitz move}. If  $\mathcal{V}=\{g_1, \ldots,
g_r; \; h_1, \ldots \, h_{2g'} \}$ is a generating vector of $G$
with respect to $\Gamma$, by definition of $\textrm{Aut}^+(\Gamma)$ any Hurwitz
move sends $g_i$ to some conjugated of $g_j$, where $o(g_i)=o(g_j)$.
In particular, if $G$ is abelian then the Hurwitz moves permute the
$g_i$ having the same order. Moreover in this case the Hurwitz moves
on $\mathcal{V}$ are unambiguously defined,
 since $\textrm{Inn}(G)$ is trivial. \\
If $\Sigma_{g'}$ is a differentiable model of a  compact Riemann surface of genus $g'$ and $p_1,\ldots, p_r \in \Sigma_{g'}$, we define
\begin{equation*}   
\textrm{Mod}_{g', [r]}:=\pi_0 \; \textrm{Diff}^+(\Sigma_{g'}-\{ p_1, \ldots,
p_r \}).
\end{equation*}
Given $\Gamma:=\Gamma(g' \; | \; m^r)$, it is well known that
$\textrm{Mod}(\Gamma)$ is isomorphic to $\textrm{Mod}_{g', [r]}$ (\cite{Schn03}, Theorem
2.2.1). In the sequel of this paper we will deal with an abelian group $G$
and with few types of signature, namely $(0 \; |\; \mathbf{m})$, $(1 \; | \; m)$ and $(1 \; | \; m^2)$. So let us explicitly describe the Hurwitz moves in these cases. 
\subsection{The case $g'=0$} \label{g'=0}
For the sake of simplicity, let us suppose that all the $m_i$ are equal,
i.e. $m_1= \ldots=m_r=m$. By the result mentioned above, the mapping class group
$\textrm{Mod}(\Gamma(0 \; | \; m^r))$ can be identified with
\begin{equation*}
\textrm{Mod}_{0,[r]}:=\pi_0 \;\textrm{Diff}^+(\mathbb{P}^1 -\{p_1,
\ldots p_r \}),
\end{equation*}
which is a quotient of the Artin braid group $\mathbf{B}_r$. Let
$\sigma_i$ be the positive-oriented Dehn twist about a simple closed
curve in $\mathbb{P}^1$ containing $p_i$ and none of the other marked points. Then it is
well known (see for instance \cite{Schn03}, Section 2.3) that
$\textrm{Mod}_{0,[r]}$ is generated by  $\sigma_1, \ldots, \sigma_r$
with the following relations:
\begin{equation*}
\begin{split}
& \sigma_i \sigma_{i+1} \sigma_i=\sigma_{i+1} \sigma_i
\sigma_{i+1}\\
& \sigma_i \sigma_j = \sigma_j \sigma_i \quad \textrm{if} \; |i-j| \geq 2\\
& \sigma_{r-1} \sigma_{r-2} \cdots \sigma_1^2 \cdots \sigma_{r-2}
\sigma_{r-1}=1.
\end{split}
\end{equation*}
Now we can describe the Hurwitz moves in this case.
\begin{proposition} \label{Among=0}
Up to inner automorphisms, the action of $\emph{Mod}_{0,[r]}$ on
$\Gamma(0 \; | \; m^r)$ is given by
\begin{equation*}
\sigma_i: \left \{
\begin{array}{ll}
y_i & \lr y_{i+1} \\
y_{i+1} & \lr y_{i+1}^{-1} y_i y_{i+1} \\
y_j & \lr y_j  \quad \emph{if} \; j \neq i,\; i+1;
\end{array} \right.
\end{equation*}
\end{proposition}
\begin{proof}
See \cite{Schn03}, Proposition 2.3.5 or \cite{Ca05},
Section 4.
\end{proof}
\begin{corollary} \label{Ahurmov}
 Let $G$ be a finite abelian group and
 let $\mathcal{V}=\{g_1, \ldots, g_r\}$ be a
generating vector of $G$ with respect to $\Gamma(0\; | \; m^r)$.
Then the Hurwitz moves coincide with the group of permutations of
$\mathcal{V}$.
\end{corollary}
The general case can be carried out in a similar way (see [Br90], Proposition 2.5.) and one obtains
\begin{corollary} \label{Aurmov}
Let $G$ be a finite abelian group and
 let $\mathcal{V}=\{g_1, \ldots, g_r\}$ be a
generating vector of $G$ with respect to $\Gamma(0\; | \; m_1,
\ldots, m_r)$. Then the Hurwitz moves on $\mathcal{V}$ are generated
by the transpositions of the $g_i$ having the same order.
\end{corollary}

\subsection{The case $g'=1, \; r=1$} \label{g'=1,n=1}
\begin{figure}[ht!]
\begin{center}
\includegraphics*[totalheight=4.5cm]{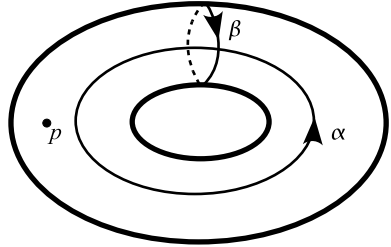}
\end{center}
\caption{Generators of $\textrm{Mod}_{1,1}$} \label{Mod11}
\end{figure}

\begin{figure}[ht!]
\begin{center}
\includegraphics*[totalheight=7cm]{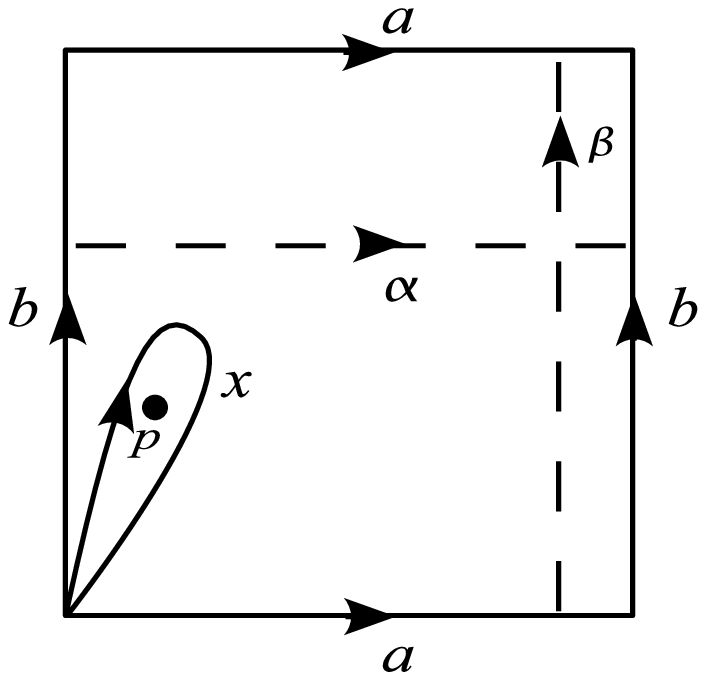}
\end{center}
\caption{} \label{Cycles11}
\end{figure}

\begin{figure}[ht!]
\begin{center}
\includegraphics*[totalheight=6.5cm]{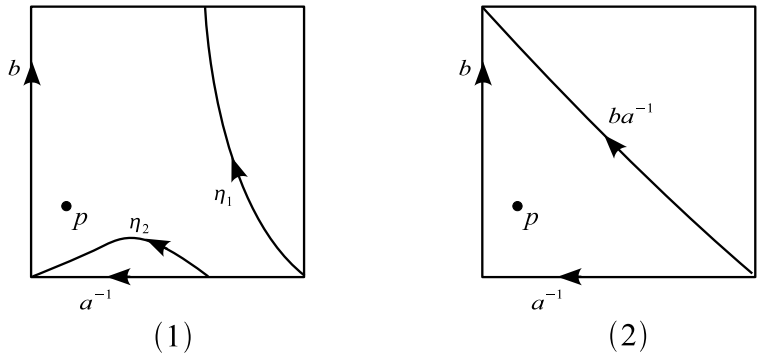}
\end{center}
\caption{$t_{\beta}(a^{-1})=\eta_1 \eta_2=ba^{-1}$}
\label{t-beta-a11}
\end{figure}

Let $\Gamma=\Gamma(1\;|\;m^1)$; then  $\textrm{Mod}(\Gamma)$ can be identified with
\begin{equation*}
\textrm{Mod}_{1,1}=\pi_0 \; \textrm{Diff}^+(\Sigma_1- \{p \}).
\end{equation*}
This group is generated by the
positively-oriented Dehn twists $t_{\alpha}, \; t_{\beta}$ about the
two simple closed curves $\alpha, \; \beta$ shown in Figure
\ref{Mod11}.
The corresponding relations are the following (see \cite{Schn03}):
\begin{equation*}
t_{\alpha}t_{\beta}t_{\alpha}=t_{\beta}t_{\alpha}t_{\beta}; \quad
(t_{\alpha}t_{\beta})^3=1.
\end{equation*}
Via the identifications
\begin{equation*}
t_{\alpha}=\left( \begin{array}{cc}
1 & 1\\
0  & 1 \\
\end{array} \right),  \quad
t_{\beta}=\left( \begin{array}{cc}
\; \; \;1 & 0\\
-1 & 1 \\
\end{array} \right)
\end{equation*}
one verifies that $\textrm{Mod}_{1,1}$ is isomorphic to
$\textrm{SL}_2(\mathbb{Z})$. The group $\Gamma(1 \;| \; m^1)$ is a
quotient of $\pi_1(\Sigma_1 -\{p \})$, in fact it has the
presentation
\begin{equation*}
\Gamma(1\; | \; m^1)=\langle a,b,x  \; | \; x^m=x[a,b]=1 \rangle.
\end{equation*}
Let us identify the torus $\Sigma_1$ with the topological space
obtained by gluing the opposite sides of a square; then the
generators $a,b,x$ of $\Gamma(1 \; | \; m^1)$ and the two loops
$\alpha, \beta$ are illustrated in Figure \ref{Cycles11}.
\begin{proposition} \label{mong1n1}
Up to inner automorphisms, the action of $\emph{Mod}_{1,1}$ on
$\Gamma(1\;| \; m^1)$ is given by
\begin{equation*}
t_{\alpha} \colon \left\{ \begin{array}{ll}
x \lr x \\
a \lr a \\
b \lr ba \\
\end{array} \right.
\quad t_{\beta} \colon \left\{ \begin{array}{ll}
x \lr x \\
a \lr ab^{-1} \\
b \lr b. \\
\end{array} \right.
\end{equation*}
\end{proposition}
\begin{proof}
It is sufficient to compute, up to inner automorphisms, the action
of $\textrm{Mod}_{1,1}$
 on $\pi_1(\Sigma_1-\{p\})$. Look at Figure \ref{Cycles11}.
Evidently, $t_{\alpha}(a)=a$ and $t_{\alpha}(x)=x$, because $\alpha$
is disjoint from both $a$ and $x$. Analogously, $t_{\beta}(b)=b, \;
t_{\beta}(x)=x$. Hence we must only compute $t_{\alpha}(b)$ and
$t_{\beta}(a)$. The pair $(\alpha, b)$ is positively oriented and
the action of $t_{\alpha}$ on $b$ is illustrated in Figure
\ref{t-alpha-b11}; so we obtain $t_{\alpha}(b)=\xi_1\xi_2$, which is
homotopic to $ba$.
\begin{figure}[ht!]
\begin{center}
\includegraphics*[totalheight=6.5cm]{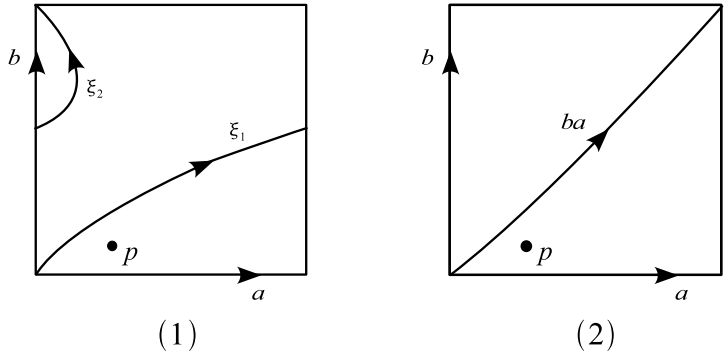}
\end{center}
\caption{$t_{\alpha}(b)=\xi_1\xi_2=ba$} \label{t-alpha-b11}
\end{figure}
Similarly, $(\beta, a^{-1})$ is a positively oriented pair and the
action of  $t_{\beta}$ on $a^{-1}$ is illustrated in Figure
\ref{t-beta-a11}; so we obtain $t_{\beta}(a^{-1})=\eta_1
\eta_2=ba^{-1}$, that is $t_{\beta}(a)=ab^{-1}$.
\end{proof}
\begin{corollary} \label{hurg1n1}
Let $G$ be a finite abelian group and let $\mathcal{W}=\{g; \;
h_1,\;h_2\}$ be a generating vector for $G$ with respect to $\Gamma(1
\;| \; m^1)$. Then the Hurwitz moves on $\mathcal{W}$ are generated
by
\begin{equation*}
 \mathbf{1} \colon \left\{\begin{array}{ll}
 g  \lr g \\
 h_1 \lr h_1 \\
 h_2 \lr h_1+h_2 \\
\end{array} \right.
\quad \mathbf{2} \colon \left\{\begin{array}{ll}
 g  \lr g \\
 h_1 \lr h_1-h_2 \\
 h_2 \lr h_2. \\
\end{array} \right.
\end{equation*}
\end{corollary}
\begin{proof}
This follows directly from Proposition \ref{mong1n1}.
\end{proof}

\subsection{The case $g'=1, \; r=2, \; m_1=m_2=m$} \label{g'=1,n=2}
\begin{figure}[ht!]
\begin{center}
\includegraphics*[totalheight=4.5cm]{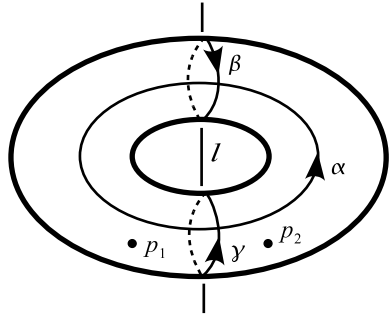}
\end{center}
\caption{Generators of $\textrm{Mod}_{1,2}$} \label{Mod12}
\end{figure}

\begin{figure}[ht!]
\begin{center}
\includegraphics*[totalheight=6.5cm]{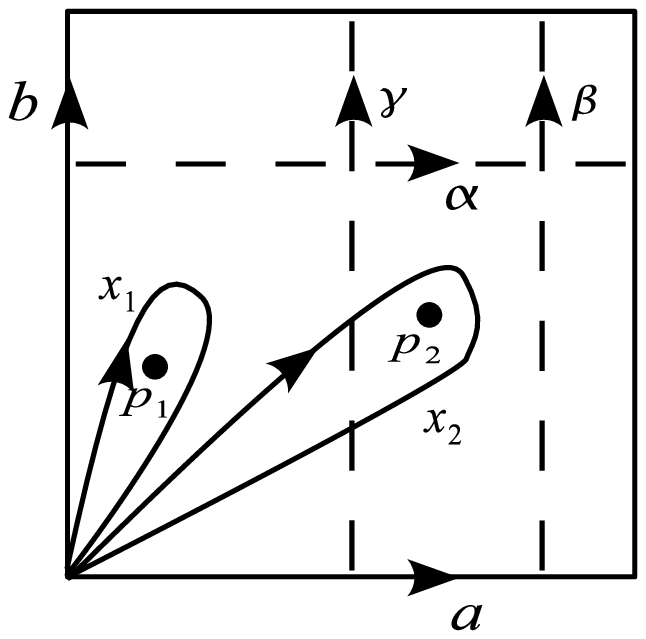}
\end{center}
\caption{} \label{Cycles12}
\end{figure}

Let $\Gamma=\Gamma(1\;|\;m^2)$; then $\textrm{Mod}(\Gamma)$ can be identified with
\begin{equation*}
\textrm{Mod}_{1,[2]}=\pi_0 \; \textrm{Diff}^+(\Sigma_1- \{p_1,p_2
\}).
\end{equation*}
This group is generated by the
positively-oriented Dehn twists $t_{\alpha}, \; t_{\beta}, \;
t_{\gamma}$ about the simple closed curves $\alpha, \; \beta, \;
\gamma$ shown in Figure \ref{Mod12}, and by the class of the
rotation $\rho$ of $\pi$ radians around the line $l$,  which
exchanges the marked points.
The relations defining $\textrm{Mod}_{1,[2]}$ are the following (see
\cite{CaMu04}):
\begin{equation*}
\begin{split}
t_{\alpha}t_{\beta}t_{\alpha}&=t_{\beta}t_{\alpha}t_{\beta}; \quad
t_{\alpha}t_{\gamma}t_{\alpha}=t_{\gamma}t_{\alpha}t_{\gamma}; \\
t_{\beta}t_{\gamma}&=t_{\gamma}t_{\beta}; \quad
(t_{\alpha}t_{\beta}t_{\gamma})^4=1;\\
t_{\alpha}\rho&=\rho t_{\alpha}; \quad t_{\beta}\rho=\rho t_{\beta};
\quad t_{\gamma}\rho=\rho t_{\gamma}.
\end{split}
\end{equation*}
The group $\Gamma(1 \;| \; m^2)$ is a quotient of
$\pi_1(\Sigma_1-\{p_1,p_2 \})$, in fact its presentation is
\begin{equation*}
\Gamma(1 \; | \; m^2)=\langle a,b,x_1,x_2
\;|\;x_1^m=x_2^m=x_1x_2[a,b]=1 \rangle;
\end{equation*}
the generators $a, \; b, \; x_1, \; x_2$ and the loops $\alpha, \;
\beta, \; \gamma$ are illustrated in Figure \ref{Cycles12}.
\begin{proposition} \label{mon1,2}
Up to inner automorphisms, the action of $\emph{Mod}_{1,[2]}$ on
$\Gamma(1 \;| \; m^2)$ is given by
\begin{equation*}
\begin{split}
& t_{\alpha} \colon \left\{ \begin{array}{ll}
  x_1 \lr x_1 \\
  x_2 \lr x_2 \\
  a \lr a \\
  b \lr ba
\end{array} \right.
\quad  \quad \quad \quad \quad \quad \quad t_{\beta} \colon \left\{
\begin{array}{ll}
  x_1 \lr x_1 \\
  x_2 \lr x_2\\
  a \lr ab^{-1} \\
  b \lr b
\end{array} \right.
\\
& t_{\gamma} \colon \left\{ \begin{array}{ll}
  x_1 \lr x_1 \\
  x_2 \lr ab^{-1}a^{-1}x_2aba^{-1}\\
  a \lr b^{-1}x_1a \\
  b \lr b
\end{array} \right.
\quad  \; \rho \colon \left\{ \begin{array}{ll}
  x_1 \lr b^{-1}a^{-1}x_2 ab \\
  x_2 \lr a^{-1}b^{-1}x_1 ba \\
  a \lr a^{-1} \\
  b \lr b^{-1}.
\end{array} \right.
\end{split}
\end{equation*}
\end{proposition}
\begin{proof}
It is sufficient to compute, up to inner automorphisms, the action
of $\textrm{Mod}_{1, [2]}$ on $\pi_1(\Sigma_1-\{p_1,p_2\})$. Look at
Figure \ref{Cycles12} and consider the action of $t_{\alpha}$. We
have $t_{\alpha}(a)=a, \; t_{\alpha}(x_1)=x_1, \;
t_{\alpha}(x_2)=x_2$ because $\alpha$ is disjoint from $a, \; x_1,
\;x_2$; moreover $t_\alpha(b)=ba$ exactly as in the proof of
Proposition \ref{mong1n1}. The computation of the action of
$t_{\beta}$ is similar. Next, let us consider the action of
$t_{\gamma}$. The curve $\gamma$ is disjoint from both $b$ and
$x_1$, then $t_{\gamma}(b)=b, \; t_{\gamma}(x_1)=x_1$. Moreover,
since $(\gamma, a^{-1})$ is a positively-oriented pair, the action
of $t_{\gamma}$ on $a^{-1}$ is as in Figure \ref{t-gamma-a12}; this
gives $t_{\gamma}(a^{-1})=\xi_1\xi_2=a^{-1}x_1^{-1}b$, hence
$t_{\gamma}(a)=b^{-1}x_1a$.
\begin{figure}[ht!]
\begin{center}
\includegraphics*[totalheight=6.5cm]{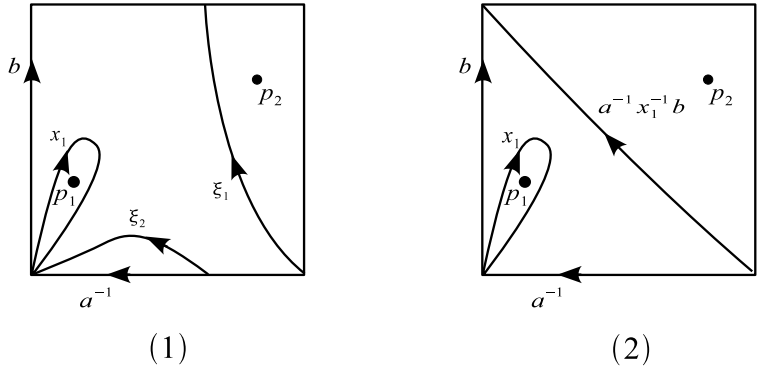}
\end{center}
\caption{$t_{\gamma}(a^{-1})= \xi_1 \xi_2=a^{-1}x_1^{-1}b$}
\label{t-gamma-a12}
\end{figure}
Using the computations above and the relation $x_1x_2[a,b]=1$ we
obtain $x_1t_{\gamma}(x_2)[b^{-1}x_1a, \;b]=1$, hence
\begin{equation*}
t_{\gamma}(x_2)=ab^{-1}a^{-1}x_1^{-1}b=ab^{-1}a^{-1}x_2aba^{-1}.
\end{equation*}
Finally, let us consider the action of $\rho$, which is illustrated
in Figure \ref{Action-rho}.
\begin{figure}[ht!]
\begin{center}
\includegraphics*[totalheight=6.5cm]{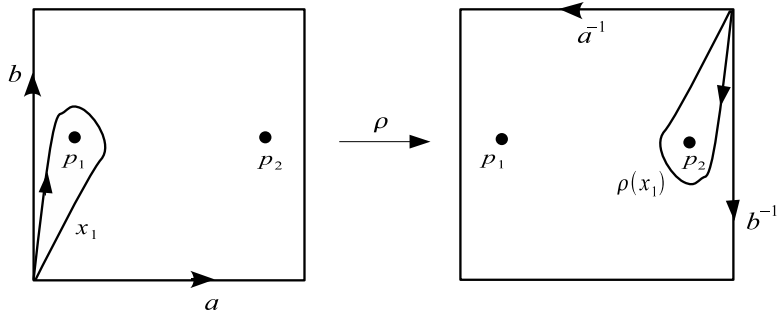}
\end{center}
\caption{Action of $\rho$} \label{Action-rho}
\end{figure}
Evidently, $\rho(a)=a^{-1}$ and $\rho(b)=b^{-1}$. Moreover, the
picture
 also shows that there is the relation
\begin{equation*}
ba(\rho(x_1))^{-1}b^{-1}a^{-1}=x_1.
\end{equation*}
Then
\begin{equation*}
\begin{split}
\rho(x_1)&=b^{-1}a^{-1}x_1^{-1}ba=b^{-1}a^{-1}(x_2aba^{-1}b^{-1})ba \\
&=b^{-1}a^{-1}x_2ab.
\end{split}
\end{equation*}
Finally, using $x_1x_2[a,b]=1$, we can write
\begin{equation*}
\begin{split}
\rho(x_2)&=(\rho(x_1))^{-1}(\rho([a,b]))^{-1} \\
&=(b^{-1}a^{-1}x_2ab)^{-1}(a^{-1}b^{-1}ab)^{-1} \\
&=b^{-1}a^{-1}x_2^{-1}ba=a^{-1}b^{-1}x_1ba.
\end{split}
\end{equation*}
This completes the proof.
\end{proof}
\begin{corollary} \label{hurmov1,2}
Let $G$ be a finite, abelian group and let $\mathcal{W}=\{g_1,\;g_2; \;
\; h_1,\;h_2 \}$ be a generating vector for $G$ with respect to
$\Gamma(1 \; | \; m^2)$. Then the Hurwitz moves on $\mathcal{W}$ are
generated by
\begin{equation*}
\begin{split}
 &  \mathbf{1} \colon \left\{ \begin{array}{ll}
  g_1 \lr g_1 \\
  g_2 \lr g_2 \\
  h_1 \lr h_1 \\
  h_2 \lr h_1+h_2
\end{array} \right.
\quad \quad \quad \mathbf{2} \colon \left\{ \begin{array}{ll}
 g_1 \lr g_1 \\
 g_2 \lr g_2 \\
 h_1 \lr h_1-h_2 \\
 h_2 \lr h_2
\end{array} \right.
\\
 & \mathbf{3} \colon \left\{ \begin{array}{ll}
  g_1 \lr g_1 \\
  g_2 \lr g_2 \\
  h_1 \lr h_1-h_2+g_1 \\
  h_2 \lr h_2
\end{array} \right.
\quad \mathbf{4} \colon \left\{ \begin{array}{ll}
  g_1 \lr g_2 \\
  g_2 \lr g_1 \\
  h_1 \lr -h_1 \\
  h_2 \lr -h_2.
\end{array} \right.
\end{split}
\end{equation*}
\end{corollary}

\section{Surfaces of general type with $p_g=q=1$ isogenous to a
 product} \label{surfisog}

\begin{definition} \label{iso-prod}
 A surface $S$ of general type is said to be \emph{isogenous to a product} if there exist two smooth curves $C, \; F$ and a finite group $G$, acting freely on
their product, so that $S=(C \times F)/G$. 
\end{definition}
We have two cases: the
\emph{mixed} case, where the action of $G$ exchanges the two factors
(and then $C, \; F$ are isomorphic) and the \emph{unmixed} case,
where $G$ acts diagonally. If $S$ is isogenous to a product, there exists a unique realization $S=(C \times F) /G$ such that the genera $g(C),
\; g(F)$ are minimal (\cite{Ca00}, Proposition 3.13). Our aim is to solve the following
\begin{prob}
Classify the surfaces of general type $S=(C \times F)/G$ with
$p_g=q=1$, isogenous to an \emph{unmixed} product, assuming that the group
$G$ is \emph{abelian}. Describe the corresponding irreducible
 components of the moduli space.
\end{prob}
Notice that, when $S$ is of unmixed type, the group
 $G$ acts  separately on $C$ and $F$ and the two projections
  $ \pi_C \colon C
\times F \lr C, \; \pi_F \colon C \times F  \lr F$ induce two
isotrivial fibrations $\alpha \colon S \lr C/G, \; \beta \colon S
\lr F/G$
 whose smooth fibres are isomorphic to $F$ and $C$, respectively. Moreover,
 we will always consider the minimal realization of $S$, so that the action
 of $G$ will be faithful on both factors. 
\begin{proposition} \label{basic facts}
Let $S=(C \times F) /G$ be a surface of general type with $p_g=q=1$,
isogenous to an unmixed product. Then $S$ is minimal and moreover
\begin{itemize}
\item[$(i)$] $K_S^2=8$;
\item[$(ii)$] $|G|=(g(C)-1)(g(F)-1)$;
\item[$(iii)$]$F/G \cong \mathbb{P}^1$ and $C/G \cong E$, where $E$ is
 an elliptic curve isomorphic to the Albanese variety of $S$;
\item[$(iv)$] $S$ contains no pencils of genus $2$ curves;
\item[$(v)$] $g(C) \geq 3, \; \; g(F) \geq 3  \quad (\textrm{hence $(ii)$ implies }|G| \geq 4)$.
\end{itemize}
\end{proposition}
\begin{proof}
Since the projection $C \times F \lr S$ is \'{e}tale, the pullback
of any $(-1)-$curve of $S$ would give rise to a (disjoint) union of
$(-1)-$curves in $C \times F$, and this is impossible; then $S$ is
minimal. Now let us prove $(i)-(v)$.
\begin{itemize}
\item[$(i)$] Since $S$ is isogenous to an unmixed product we have
$K_S^2=2c_2(S)$ (\cite{Se90}, Proposition 3.5). Noether's formula
 gives $K_S^2+c_2(S)=12$, so it follows  $K^2_S=8$.
\item[$(ii)$] We have
\begin{equation*}
p_g(C \times F)=g(C) \cdot g(F) \quad \textrm{and} \quad
q(C \times F)=g(C)+g(F),
\end{equation*}
see \cite{Be2}, III.22. Since $C \times F
\longrightarrow S$ is an {\'e}tale covering, we obtain $|G|=
\chi(\mO_{C \times F}) / \chi(\mO_S)=(g(C)-1)(g(F)-1)$.
\item[$(iii)$] We have $q(S)=g(C/G)+g(F/G)$, then we may assume
\begin{equation*}
g(C/G)=1, \;\; g(F/G)=0.
\end{equation*}
Setting $E=C/G$ it follows that $\alpha \colon S \longrightarrow E$
is a connected fibration with elliptic base, then it coincides with
the Albanese morphism of $S$.
\item[$(iv)$] Let $S$ be a minimal surface of general type with
$p_g=q=1$ which contains a genus $2$ pencil. There are two cases:
\begin{itemize}
\item[$\bullet$] the pencil is rational, then either $K_S^2=2$ or
  $K_S^2=3$ (see \cite{Xi85}, p. 51);
\item[$\bullet$] the pencil is irrational, therefore it must be the Albanese
pencil, and in this case $2 \leq K_S^2 \leq 6$ (see \cite{Xi85}, p.
17).
\end{itemize}
In both cases, part $(i)$ implies that $S$  cannot be isogenous to a
product.
\item[$(v)$] This follows from part $(iv)$.
\end{itemize}
\end{proof}
The two $G-$coverings $f \colon F \lr \mathbb{P}^1$, $ h \colon C
\lr E$
 are induced by two admissible epimorphisms $\vartheta \colon \Gamma(0 \; | \; \mathbf{m}) \lr G$, $\psi \colon  \Gamma(1 \; | \; \mathbf{n}) \lr G$,
where $\mathbf{m}=(m_1, \ldots, m_r), \; \mathbf{n}=(n_1, \ldots,
n_s)$. The Riemann-Hurwitz formula gives
 \begin{equation} \label{generi1}
\begin{split}
2g(F)-2 &=|G|\bigg(-2+\sum_{i=1}^r \bigg( 1- \frac{1}{\;m_i} \bigg)
\bigg) \\
2g(C)-2 & =|G| \sum_{j=1}^s \bigg(1-\frac{1}{\;n_j} \bigg).
\end{split}
\end{equation}
\begin{proposition} \label{various-cases}
We have the following possibilities:
\begin{itemize}
\item[$(a)$] $g(F)=3, \quad \mathbf{n}=(2^2)$
\item[$(b)$] $g(F)=4, \quad \mathbf{n}=(3^1)$
\item[$(c)$] $g(F)=5, \quad \mathbf{n}=(2^1)$.
\end{itemize}
Moreover, if $G$ is abelian only case $(a)$ may occur.
\end{proposition}
\begin{proof}
Using (\ref{generi1}) and part $(ii)$ of Proposition \ref{basic
facts} we obtain
\begin{equation} \label{TTT}
2=(g(F)-1)\sum_{j=1}^s \bigg(1-\frac{1}{\;n_j} \bigg).
\end{equation}
Since $g(F) \geq 3$ and the sum in the right-hand side of
(\ref{TTT}) is $\geq \frac{1}{2}$, we have $g(F) \leq 5$. If
$g(F)=3$ then $\sum_{j=1}^s \left(1-\frac{1}{\;n_j} \right)=1$, so
$\mathbf{n}=(2^2)$; if $g(F)=4$ then $\sum_{j=1}^s
\left(1-\frac{1}{\;n_j} \right)=\frac{2}{3}$, so $\mathbf{n}=(3^1)$;
 if $g(F)=5$ then  $\sum_{j=1}^s \left(1-\frac{1}{\;n_j} \right)=\frac{1}{2}$,
 so $\mathbf{n}=(2^1)$. Finally, if $G$ is abelian then $s \geq 2$ (Proposition
 \ref{ab-no1}), so only case $(a)$ may occur.
 \end{proof}
\begin{remark}
If $G$ is not abelian, then all possibilities $(a)$, $(b)$, $(c)$
actually occur. Examples are given in Section \ref{newexamples}.
\end{remark}
\begin{corollary} \label{ord G}
Suppose that $G$ is abelian. Then
\begin{itemize}
\item the Albanese fibration of $S$ is an isotrivial genus $3$ pencil with
two double fibres;
\item $|G|=2(g(C)-1)$.
\end{itemize}
\end{corollary}

\section{Abelian case: the building data} \label{building}

In the sequel we will assume that $G$ is abelian. By Proposition
\ref{various-cases} the covering $h \colon C \lr E$ is induced by an
admissible epimorphism $\psi \colon \Gamma(1 \; | \; 2^2) \lr G$. If
$\mathcal{W}=\{ g_1,\; g_2; \; h_1,\; h_2 \}$ is the corresponding
generating vector, we have $2g_1=2g_2=g_1+g_2=0$, hence $g_1=g_2$.
For the sake of simplicity we set $g_1=g_2=g$ and we denote the
generating vector by $\{g ; \; h_1,\; h_2 \}$. Note that $\langle g
\rangle$ is the only non trivial stabilizer of the action of $G$ on
$C$. Analogously, if $\mathcal{V}:=\{ g_1, \ldots g_r\}$ is any generating
vector of $G$ with respect to $\Gamma(0 \; | \; \mathbf{m})$, the
cyclic subgroups $\langle g_1 \rangle, \ldots, \langle g_r \rangle$
are the only non trivial stabilizers of the action of $G$ on $F$.
Then the diagonal action of $G$ on $C \times F$ is free if and only
if
\begin{equation} \label{stab free}
\left( \bigcup_{i=1}^r \langle g_i \rangle \right) \cap \langle g
\rangle= \{0 \}.
\end{equation}
Using the results contained in the previous section, we obtain
\begin{proposition} \label{structureresult1}
Suppose that we have the following data:
\begin{itemize}
\item a finite abelian group $G$;
\item two admissible epimorphisms
\begin{equation*}
\begin{split}
& \vartheta \colon \Gamma(0 \; | \; \mathbf{m}) \lr G,  \quad \mathbf{m}=(m_1, \ldots, m_r)\\
& \psi \colon \Gamma(1 \; | \; 2^2) \lr G
\end{split}
\end{equation*}
with corresponding generating vectors  $\mathcal{V}=\{g_1, \ldots,
g_r \}, \; \mathcal{W}=\{g; \; h_1,\;h_2 \}$.
\end{itemize}
Let
\begin{equation*}
\begin{split}
& f \colon F \lr \mathbb{P}^1=F/G \\
& h \colon C \lr E= C/G
\end{split}
\end{equation*}
be the  $G-$coverings induced by $\vartheta$ and $\psi$ and let
$g(C)$, $g(F)$ be the genera of $C$ and $F$, which are related on
$G$ and $\mathbf{m}$ by $(\ref{generi1})$. Assume moreover that
\begin{itemize}
\item $g(C) \geq 3, \; \; g(F)=3$; \item $|G|=2(g(C)-1)$; \item
condition $(\ref{stab free})$ is satisfied.
\end{itemize}
Then the diagonal action of $G$ on $C \times F$ is free and the
quotient $S=(C \times F)/G$ is a minimal surface of general type
with $p_g=q=1$. Conversely, any surface of general type with
$p_g=q=1$, isogenous to an unmixed product with $G$ abelian, arises
in this way.
\end{proposition}
We will call the $4-$ple $(G, \; \mathbf{m}, \; \vartheta, \; \psi)$ the
\emph{building data} of  $S$.

\begin{corollary} \label{nociclico}
Let $S=(C \times F) /G$ be a surface of general type with $p_g=q=1$,
isogenous to an unmixed product. Then the group $G$ cannot be
cyclic.
\end{corollary}
\begin{proof}
By contradiction suppose that $G$ is cyclic; then
 $G=\mathbb{Z}_{2m}$ for some integer $m$. Let
$\{g_1, \ldots, g_r \}$ and $\{g; \; h_1, \; h_2\}$ be generating
vectors of $G$ as in Proposition \ref{structureresult1}. The group
$G$ contain exactly one subgroup of order $2$, namely $\langle g
\rangle$.
 On the other hand, since $\langle g_1, \ldots ,g_r \rangle=G$, we have
l.c.m.$(m_1,\ldots,m_r)=2m$; hence $2$ divides some of the $m_i$,
say $m_1$. This implies $g \in \langle g_1 \rangle$, which
 violates condition (\ref{stab free}).
\end{proof}

Let $\mathfrak{M}_{a,b}$ be
the moduli space of smooth minimal surfaces of general type with
$\chi(\mO_S)=a, \; K_S^2=b$; by a result of Gieseker, we know that
$\mathfrak{M}_{a,b}$ is a quasiprojective variety for all $a,b \in
\mathbb{N}$ (see \cite{Gie77}). Obviously, our surfaces are contained in $\mathfrak{M}_{1,8}$
and we want to describe their locus there.
\begin{proposition} \label{moduli-comp}
For fixed $G$ and
  $\mathbf{m}$, denote by
$\mathfrak{M}(G, \; \mathbf{m})$
 the moduli space of surfaces with $p_g=q=1$ described by
 the remaining building data $(\vartheta,\; \psi)$. Then
$\mathfrak{M}(G, \; \mathbf{m})$ consists of a finite number of
 irreducible components of $\mathfrak{M}_{1,8}$, all of dimension $r-1$.
\end{proposition}
\begin{proof}
The fact that $\mathfrak{M}(G, \; \mathbf{m})$ consists of finitely
many irreducible components of $\mathfrak{M}_{1,8}$ follows by general results
of Catanese on surfaces isogenous to a product
 (see \cite{Ca00}). The dimension of each component is $r-1$
because we take $r$ points on $\mathbb{P}^1$ modulo projective
equivalence and $2$ points on $E$ modulo projective equivalence.
\end{proof}
Let us define
\begin{equation*}
\begin{split}
\Phi(G,\; \mathbf{m})&:= \textrm{Epi}(\Pi_3, \; \Gamma(0 \;|\;
\mathbf{m}), \; G) \times  \textrm{Epi}(\Pi_{g(C)}, \; \Gamma(1
\;|\;
2^2), \; G); \\
\mathfrak{G}&:=\textrm{Aut}(G) \times \textrm{Mod}(\Gamma(0 \;|\;
\mathbf{m})) \times \textrm{Mod}(\Gamma(1 \;|\; 2^2)).
\end{split}
\end{equation*}
The group $\mathfrak{G}$ naturally acts on the set $\Phi(G,\;
\mathbf{m})$ in the following way:
\begin{equation*}
(\lambda, \; \eta_0, \; \eta_1) \cdot (\vartheta, \; \psi):=(\lambda
\circ \vartheta \circ \eta_0, \; \lambda \circ \psi \circ \eta_1).
\end{equation*}
\begin{proposition} \label{structuremoduli}
Let $S_1, \;S_2$ be two surfaces defined by building data
$(\vartheta_1, \; \psi_1), \; (\vartheta_2, \; \psi_2) \in \Phi(G,
\; \mathbf{m})$. Then $S_1$ and $S_2$ belong to the same connected
component of $\mathfrak{M}(G, \; \mathbf{m})$ if and only if
$(\vartheta_1, \; \psi_1)$ and $(\vartheta_2, \; \psi_2)$ are in the
same $\mathfrak{G}-$class.
\end{proposition}
\begin{proof}
We can use the same argument of \cite{BaCa03}, Theorem 1.3. In fact,
Proposition \ref{same-class} allows us to substitute the pair of
braid group actions considered in that paper with the two actions of
$\textrm{Mod}(\Gamma(0 \;|\; \mathbf{m}))$ and
$\textrm{Mod}(\Gamma(1 \;|\; 2^2))$.
\end{proof}
Now let $\mathfrak{B}(G, \; \mathbf{m})$ be the set of
 pairs of generating vectors $(\mathcal{V}, \; \mathcal{W})$
  such that the the hypotheses of Proposition \ref{structureresult1} are
  satisfied (in particular (\ref{stab free}) must hold).
 Let us denote by $\mathfrak{R}$ the equivalence relation on
$\mathfrak{B}(G, \; \mathbf{m})$ generated by
\begin{itemize}
\item Hurwitz moves on $\mathcal{V}$;
\item Hurwitz moves on $\mathcal{W}$;
\item simultaneous conjugation of $\mathcal{V}$ and $\mathcal{W}$ by an element
of $\lambda \in \textrm{Aut}(G)$, i.e. we let $(\mathcal{V}, \;
\mathcal{W})$ be equivalent to $(\lambda(\mathcal{V}), \;
\lambda(\mathcal{W}))$.
\end{itemize}

\begin{proposition} \label{structuremoduli1}
The number of irreducible components in $\mathfrak{M}(G, \;
\mathbf{m})$ equals the number of $\mathfrak{R}-$classes in
$\mathfrak{B}(G, \; \mathbf{m})$.
\end{proposition}
\begin{proof}
Immediate consequence of Proposition \ref{structuremoduli}.
\end{proof}

\section{Abelian case: the classification}  \label{classification}

We have the following result.

\begin{theorem} \label{isog prod}
If  the group $G$ is abelian, then there exist exactly four families
of surfaces of general type with $p_g=q=1$, isogenous to an unmixed
product. In any case $g(F)=3$, whereas the possibilities for $g(C)$
and $G$ are the following:
\begin{itemize}
\item[$I$.] $g(C)=3, \quad G=(\mathbb{Z}_2)^2$;
\item[$II$.] $g(C)=5, \quad G=(\mathbb{Z}_2)^3$;
\item[$III$.] $g(C)=5, \quad G=\mathbb{Z}_2 \times \mathbb{Z}_4$;
\item[$IV$.] $g(C)=9, \quad G= \mathbb{Z}_2 \times \mathbb{Z}_8$.
\end{itemize}
\end{theorem}
Surfaces of type $I$ already appear in  \cite{Pol06}, whereas those
of type $II, \; III, \; IV$ provide new examples of minimal surfaces
of general type with $p_g=q=1, \; K^2=8$. The remainder of this
section deals with the proof of Theorem \ref{isog prod}. Let $S$ be
defined by building data $(G, \; \mathbf{m}, \; \vartheta, \; \psi)$
as in Proposition \ref{structureresult1}. Using (\ref{generi1}) and
Corollary \ref{ord G} we obtain
\begin{equation} \label{diofantina}
2=(g(C)-1) \left( -2+\sum_{i=1}^r \left( 1-\frac{1}{\;m_i} \right)
\right).
\end{equation}
\begin{proposition} \label{r=6}
We have $3 \leq r \leq 6$. Moreover, if $r=6$ the only possibility is
\begin{equation*}
\mathbf{m}=(2^6), \quad G=(\mathbb{Z}_2)^2.
\end{equation*}
\end{proposition}
\begin{proof}
Using equation (\ref{diofantina}) we can write
\begin{equation} \label{stima r6 1}
 (g(C)-1) \bigg( \frac{r}{2}-2 \bigg) \leq 2 < (g(C)-1)(r -2).
 \end{equation}
Part $(v)$ of Proposition \ref{basic facts} yields $g(C)-1 \geq 2$,
 hence (\ref{stima r6 1}) implies $3 \leq r
 \leq 6$. Moreover we have $r=6$ if and only if
 $\mathbf{m}=(2^6)$, and in this case $|G|=4$.
 Then Corollary $\ref{nociclico}$ implies $G=(\mathbb{Z}_2)^2$.
\end{proof}
\begin{proposition} \label{r=5}
If $r=5$ the only possibility is
 \begin{equation*}
 \mathbf{m}=(2^5), \quad G=(\mathbb{Z}_2)^3.
 \end{equation*}
 \end{proposition}
\begin{proof}
If $r=5$ formula (\ref{diofantina}) gives
\begin{equation} \label{stima r51}
2= (g(C)-1) \bigg( 3-\sum_{i=1}^{5} \frac{1}{\;m_i} \bigg)
 \geq  (g(C)-1) \bigg( 3-\frac{5}{\;m_1} \bigg),
\end{equation}
hence
\begin{equation*}
g(C)-1 \leq \frac{2m_1}{3m_1-5} \; \cdot
\end{equation*}
If $m_1 \geq 3$ then $g(C) \leq 2$, a contradiction. Then $m_1=2$
and
  $g(C)-1 \leq 4$, hence  $|G| \leq 8$  (Corollary
 \ref{ord G}) with equality if and only if $\mathbf{m}=(2^5)$.
If $\mathbf{m} \neq (2^5)$ then $G$ would be a noncyclic group of
order smaller
 than $8$ which contains some element of order greater than $2$, a contradiction.
Therefore $\mathbf{m}=(2^5)$ is actually the only possibility.
 Then $G$ is an abelian group of order
$8$ generated by elements of order $2$, hence $G=(\mathbb{Z}_2)^ 3$.
\end{proof}
\begin{proposition} \label{r=4} If $r=4$ the only possibility
is
\begin{equation*}
\mathbf{m}=(2^2,4^2), \quad G=\mathbb{Z}_2 \times \mathbb{Z}_4.
\end{equation*}
\end{proposition}
\begin{proof}
If $r=4$ formula (\ref{diofantina}) gives
\begin{equation} \label{stima r41}
2=(g(C)-1) \bigg( 2-\sum_{i=1}^4 \frac{1}{\;m_i} \bigg) \geq
 (g(C)-1) \bigg(2-\frac{4}{\;m_1} \bigg),
\end{equation}
hence
\begin{equation*}
g(C)-1 \leq \frac{m_1}{m_1-2} \; \cdot
\end{equation*}
If $m_1 \geq 3$ then $g(C)-1 \leq 3$, which implies $|G| \leq 6$.
Since $G$ is not cyclic the only possibility would be
$G=\mathbb{Z}_2 \times \mathbb{Z}_2$, which contains no elements of
order $\geq 3$, a contradiction. It follows $\mathbf{m}=(2,
m_2,m_3,m_4)$. Applying again formula (\ref{diofantina}) we get
\begin{equation} \label{stima r42}
2=(g(C)-1) \bigg( \frac{3}{2} - \sum_{i=2}^4 \frac{1}{\;m_i} \bigg)
\geq (g(C)-1) \bigg( \frac{3}{2} - \frac{3}{\;m_2} \bigg),
\end{equation}
hence
\begin{equation*}
g(C)-1 \leq \frac{4m_2}{3m_2-6} \; \cdot
\end{equation*}
 If $m_2 \geq 3$ then $|G| \leq 8$, and equality
 occurs if and only if $\mathbf{m}=(2,3^3)$. Since a group of order
 $8$  contains no elements of order $3$, we must have $|G|<8$
and the only possibility would be $G=\mathbb{Z}_2 \times
 \mathbb{Z}_2$, which contradicts $m_2 \geq 3$. Then
 $m_2=2$, so that  $\mathbf{m}=(2^2,m_3,m_4)$. Now suppose $m_3=2$ and consider
 the generating vector $\mathcal{V}=\{g_1,g_2,g_3,g_4\}$; since
 $g_1+g_2+g_3+g_4=0$ it follows $m_4=2$,
 hence $\mathbf{m}=(2^4)$ which violates (\ref{diofantina}).
Therefore $m_3 \geq 3$ and an easy computation
  using (\ref{diofantina}) shows that there are only the following
  possibilities:
\begin{itemize}
\item[$(i)$] $\mathbf{m}=(2^2,3^2), \; \; \; \; \quad |G|=12$;
\item[$(ii)$] $\mathbf{m}=(2^2,3,6), \; \;  \quad |G|=8$;
\item[$(iii)$] $\mathbf{m}=(2^2,4^2), \; \; \; \; \quad |G|=8$;
\item[$(iv)$] $\mathbf{m}=(2^2,4,12), \quad |G|=6$;
\item[$(v)$] $\mathbf{m}=(2^2,6^2), \; \; \; \; \quad |G|=6$.
\end{itemize}
 In case $(i)$, equality $g_1+g_2+g_3+g_4=0$
yields $g_3+g_4=g_1+g_2$, hence $g_1+g_2=g_3+g_4=0$. This implies
$G= \langle g_1,g_3 \rangle=\mathbb{Z}_2 \times \mathbb{Z}_3$, which
is a contradiction. Case $(ii)$ must be excluded since a group of
order $8$ contains no elements of order $3$. Finally, cases $(iv)$
and $(v)$ must be excluded since an abelian group of order $6$ is
cyclic. Then the only possibility is $(iii)$, so
 $|G|$ is a noncyclic group of order $8$ which contains some
elements of order $4$. It follows $G=\mathbb{Z}_2 \times
\mathbb{Z}_4$.
\end{proof}
\begin{proposition} \label{r=3}
If $r=3$ the only possibility is
\begin{equation*}
 \mathbf{m}=(2,8^2), \quad G=\mathbb{Z}_2 \times \mathbb{Z}_8.
\end{equation*}
\end{proposition}
\begin{proof}
If $r=3$ formula (\ref{diofantina}) gives
\begin{equation} \label{stima r3 1}
2=(g(C)-1) \bigg( 1-\sum_{i=1}^3 \frac{1}{\;m_i} \bigg) \geq
(g(C)-1) \bigg( 1-\frac{3}{\;m_1} \bigg),
\end{equation}
hence
\begin{equation*}
g(C)-1 \leq \frac{2m_1}{m_1-3} \; \cdot
\end{equation*}
Let us consider now the generating vector $\mathcal{V}=\{g_1,g_2,g_3\}$. \\
$\mathbf{Case\; 1.\;}$ Suppose $m_1 \geq 4$. Then  $g(C)-1 \leq 8$,
so $|G| \leq 16$ and equality holds if and only if
$\mathbf{m}=(4^3)$. In this case the abelian group $G$ is generated
by two elements of order $4$, thus $G= \mathbb{Z}_4 \times
\mathbb{Z}_4$. Without loss of generality, we may suppose $g_1=e_1,
\; g_2=e_2, \; g_3=-e_1-e_2$, where $e_1=(1,0)$ and $e_2=(0,1)$.
Therefore $\langle g_1 \rangle \cup
 \langle g_2 \rangle \cup \langle g_3 \rangle$ contains all the
elements of order $2$ in $G$, and condition (\ref{stab free}) cannot
be satisfied; hence $\mathbf{m}=(4^3)$ must be excluded. It follows
$|G| <16$; since
 $G$ is not cyclic and $|G|$ is even, we are left with few possibilities.
\begin{itemize}
\item $G= \mathbb{Z}_2 \times
\mathbb{Z}_6$. This gives $m_1=6$, hence $g(C)-1 \leq 4$ and
 $|G| \leq 8$, a contradiction.
\item $G=\mathbb{Z}_2 \times \mathbb{Z}_4$. This implies that the highest order of
an element of $G$ is $4$, so $\mathbf{m}=(4^3)$, again a
contradiction.
\item $G=\mathbb{Z}_2 \times \mathbb{Z}_2$. Impossible because $m_1 \geq
4$.
\end{itemize}
Therefore  $m_1 \geq 4$ does not occur. \\
$\mathbf{Case \;2.\;}$ Suppose $m_1=3$. Since $G=\langle g_1, g_2
\rangle$ and $G$ is not cyclic, it follows $G=\mathbb{Z}_3 \times
\mathbb{Z}_{m_2}$ with $3|m_2$. Moreover $g_1+g_2+g_3=0$ implies
  $m_3=o(g_1+g_2)=m_2$. Set $m_2=m_3=m$; by using
(\ref{diofantina}) we obtain
\begin{equation} \label{stima r3 2}
2=(g(C)-1) \bigg(1 - \frac{1}{3} - \frac{2}{m} \bigg)=(g(C)-1)
\frac{2m-6}{3m}.
\end{equation}
On the other hand
\begin{equation}\label{stima r3 3}
g(C)-1= \frac{1}{2} |G| = \frac{3m}{2}.
\end{equation}
From (\ref{stima r3 2}) and (\ref{stima r3 3}) it follows $m=5$, a
contradiction.
Then $m_1=3$ cannot occur. \\
$\mathbf{Case \;3.\;}$ Suppose $m_1=2$. Exactly as before we get
$\mathbf{m}=(2,m,m)$ and $G = \mathbb{Z}_2 \times \mathbb{Z}_m$,
with $2|m$. Therefore we have
\begin{equation} \label{stima r3 4}
2=(g(C)-1) \bigg(1 - \frac{1}{2} - \frac{2}{m} \bigg)=(g(C)-1)
\frac{m-4}{2m}
\end{equation}
and
\begin{equation} \label{stima r3 5}
g(C)-1= \frac{1}{2} |G|=m.
\end{equation}
Relations (\ref{stima r3 4}) and (\ref{stima r3 5}) imply $m=8$,
hence $\mathbf{m}=(2,8^2)$ and $G=\mathbb{Z}_2 \times \mathbb{Z}_8$.
\end{proof}
This completes the proof of Theorem \ref{isog prod}.

\section{Abelian case: the moduli spaces}  \label{secmodulispaces}
Now we provide the effective construction of surfaces of type $I, \;
II, \; III, \; IV$ and the description of their moduli spaces
$\mathfrak{M}_{I},\; \mathfrak{M}_{II}, \; \mathfrak{M}_{III}, \;
\mathfrak{M}_{IV}$.
\begin{theorem} \label{modulispaces}
The moduli spaces $\mathfrak{M}_I, \; \mathfrak{M}_{II}, \;
\mathfrak{M}_{IV}$ are irreducible of dimension $5,4,2$,
respectively.
 The moduli space $\mathfrak{M}_{III}$ is the disjoint union of two irreducible
components $\mathfrak{M}_{III}^{(1)}, \;\mathfrak{M}_{III}^{(2)}$,
both of dimension $3$.
\end{theorem}
The rest of this section deals with the proof of Theorem
\ref{modulispaces}. By Proposition \ref{moduli-comp} we only have to
compute the number of irreducible components in each case; this will be done by using Proposition \ref{structuremoduli1}. Let $(\mathcal{V}, \; \mathcal{W}) \in \mathfrak{B}(G, \; \mathbf{m})$; then
 the Hurwitz moves on $\mathcal{V}=\{g_1, \ldots, g_r\}$
 are generated by the transpositions of elements having the same order (Corollary
 \ref{Aurmov}), whereas the Hurwitz moves on $\mathcal{W}=\{g; \; h_1, \; h_2 \}$ are generated by
\begin{equation*}
\begin{split}
\{g; \; h_1, \; h_2 \} & \stackrel{\mathbf{1}}{\lr} \{g; \; h_1, \;
h_1+h_2 \}, \quad \quad \, \, \, \{g; \; h_1, \; h_2 \}
\stackrel{\mathbf{2}}{\lr} \{g; \; h_1-h_2, \; h_2\},
\\
\{g; \; h_1, \; h_2 \} & \stackrel{\mathbf{3}}{\lr} \{g; \; g+
h_1-h_2, \; h_2 \}, \quad \{g; \; h_1, \; h_2 \}
\stackrel{\mathbf{4}}{\lr} \{g; \; -h_1, \; -h_2 \}
\end{split}
\end{equation*}
(see Corollary \ref{hurmov1,2}). Moreover, we will often use the Hurwitz
move obtained by successively applying $\mathbf{1}$, $\mathbf{2}$,
$\mathbf{1}$,
 and that for the sake of shortness will be denoted by $\mathbf{5}$:
\begin{equation*}
\{g; \; h_1, \; h_2 \} \stackrel{\mathbf{5}}{\lr} \{g; \; -h_2, \;
h_1 \}.
\end{equation*}

\subsection{Surfaces of type $I$}
$G=(\mathbb{Z}_2)^2, \quad \mathbf{m}=(2^6), \quad g(C)=3$. \\  Let
$\{e_1, \; e_2\}$ be the canonical basis of $G$ and consider the
generating vector  $\mathcal{W}:= \{g; \; h_1,\; h_2 \}$. Up to 
Hurwitz move $\mathbf{5}$ we may assume $\langle g, \; h_1
\rangle =G$. Modulo automorphisms of $G$ we have $g=e_1$ and
$h_1=e_2$, so there are four possibilities:
\begin{equation*}
\begin{split}
\mathcal{W}_1&=\{e_1; \;e_2, \; 0 \}, \quad \; \mathcal{W}_2=\{e_1;
\;e_2,\; e_1 \}, \\
\mathcal{W}_3&=\{e_1; \;e_2, \; e_2 \}, \quad  \mathcal{W}_4=\{e_1;
\;e_2, \; e_1+e_2 \}.
\end{split}
\end{equation*}
These vectors are all equivalent to $\mathcal{W}_1$ via a finite
sequence of Hurwitz moves:
\begin{equation*}
\mathcal{W}_2 \stackrel{\mathbf{1, \; 3, \;5, \; 3 }}{\lr}
\mathcal{W}_1; \quad \mathcal{W}_3 \stackrel{\mathbf{1}}{\lr}
\mathcal{W}_1; \quad \mathcal{W}_4 \stackrel{\mathbf{1,\; 1, \; 3,
\;5, \; 3 }}{\lr} \mathcal{W}_1.
\end{equation*}
Now let us consider $\mathcal{V}=\{g_1, \ldots, g_6 \}$. 
Condition (\ref{stab free}) implies $g_i \neq
e_1$, so there are two possibilities up to permutations:
\begin{equation*}
\begin{split}
\mathcal{V}_1&= \{e_2,\; e_2,\; e_2,\; e_2,\; e_1+e_2,\; e_1+e_2 \},\\
 \mathcal{V}_2&=\{e_1+e_2, \; e_1+e_2, \; e_1+e_2, \; e_1+e_2, \; e_2, \; e_2 \}.
\end{split}
\end{equation*}
The automorphism of $G$ given by $e_1 \to e_1, \; e_2 \to e_1+e_2$
sends $\mathcal{V}_1$ to $\mathcal{V}_2$ and $\mathcal{W}_1$ to
$\{e_1; \; e_1+e_2, \;0 \}$, which is equivalent to $\mathcal{W}_1$
via the Hurwitz move $\mathbf{3}$. This shows that the elements of $\mathfrak{B}(G;
\; \mathbf{m})$ are all $\mathfrak{R}-$equivalent, 
hence $\mathfrak{M}_I$ is
irreducible.

\subsection{Surfaces of type $II$ }
$G=(\mathbb{Z}_2)^3$, \quad
$\mathbf{m}=(2^5)$, \quad $g(C)=5$. \\
In this case the generating vector $\mathcal{W}:=\{g; \; h_1, \; h_2
\}$ must be a basis of $G$ as a $\mathbb{Z}_2-$vector space. Therefore
up to automorphisms we may assume $\mathcal{W}=\{e_1; \; e_2, \;
e_3\}$, where
 the $e_i$ form the canonical basis of $G$. Now let us consider the generating vector
$\mathcal{V}=\{g_1, \ldots, g_5 \}$; notice that condition (\ref{stab free}) implies
$g_i \neq e_1$. Let $\lambda_0$ be the automorphism of $G$ given
by $\lambda_0(e_1):=e_1$, $\; \lambda_0(e_2):=e_3$, $\; \lambda_0(e_3):=e_2$. It sends $\mathcal{W}$ to $\{e_1;\; e_3,\; e_2\}$, which is equivalent to $\mathcal{W}$ via the Hurwitz move $\mathbf{5}$. Since $\langle g_1, \ldots, g_5
\rangle=G$ and $g_1+ \cdots +g_5=0$, up to $\lambda_0$ and permutations
$\mathcal{V}$ must be one of the following:
\begin{equation*}
\begin{split}
\mathcal{V}_1 &= \{e_2, \; e_2, \; e_3, \; e_1+e_2, \; e_1+e_2+e_3 \}, \\
\mathcal{V}_2 &= \{e_2, \; e_2, \; e_1+e_3, \; e_1+e_2, \; e_2+e_3 \}, \\
\mathcal{V}_3 &= \{e_1+e_2+e_3, \; e_1+e_2+e_3, \; e_3, \; e_2+e_3, \; e_2
\}, \\
\mathcal{V}_4 &= \{e_1+e_2+e_3, \;e_1+e_2+e_3, \; e_1+e_2, \; e_2+e_3, \; e_1+e_3 \}, \\
\mathcal{V}_5 &= \{e_1+e_2, \; e_1+e_2, \; e_3, \; e_2, \; e_2+e_3 \}, \\
\mathcal{V}_6 &= \{e_1+e_2, \; e_1+e_2, \; e_1+e_3, \; e_2, \; e_1+e_2+e_3 \}, \\
\mathcal{V}_7 &= \{e_2+e_3, \; e_2+e_3, \; e_2, \; e_1+e_2+e_3, \; e_1+e_3 \}.
\end{split}  
\end{equation*}
Set $\mathcal{V}_i=\{ \alpha_i, \; \alpha_i, \; \beta_i, \; \gamma_i, \;
\delta_i \}$. One checks that, for every $i \in \{ 1, \ldots, 7 \}$, the element
$\lambda_i \in \textrm{Aut}(G)$ defined by $\lambda_i(e_1):=e_1$, $\;\lambda_i(e_2):=\alpha_i$, $\; \lambda_i(e_3):=\beta_i$
sends $\mathcal{V}_1$ to $\mathcal{V}_i$. To prove that $\mathfrak{M}_{II}$
is irreducible it is therefore sufficient to show that, for every $i$, the generating
vector $\lambda_i(\mathcal{W})$ is equivalent to $\mathcal{W}$ via a sequence
of Hurwitz moves. But this is a straightforward computation:
\begin{equation*}
\begin{tabular}{c|c|c}
$ $ & $\lambda_i(\mathcal{W})$ & $\lambda_i(\mathcal{W}) \to \mathcal{W}$
\\
\hline
$\lambda_1$ & $\{e_1; \; e_2, \; e_3 \}$ & $ $ \\
$\lambda_2$ & $\{e_1; \; e_2, \; e_1+e_3 \}$ & $\mathbf{5}, \; \mathbf{3},
\; \mathbf{2}, \; \mathbf{5}$ \\
$\lambda_3$ & $\{e_1; \; e_1+e_2+e_3, \; e_3 \}$ & $\mathbf{3}$ \\
$\lambda_4$ & $\{e_1; \; e_1+e_2+e_3, \; e_1+e_2 \}$ & $\mathbf{2}, \; \mathbf{5},
\; \mathbf{3}, \; \mathbf{2}$ \\
$\lambda_5$ & $\{e_1; \; e_1+e_2, \; e_3 \}$ & $\mathbf{3}, \; \mathbf{2}$ \\
$\lambda_6$ & $\{e_1; \; e_1+e_2, \; e_1+e_3 \}$ & $\mathbf{1}, \; \mathbf{3}, \; \mathbf{5}, \mathbf{2}$ \\
$\lambda_7$ & $\{e_1; \; e_2+e_3, \; e_2 \}$ & $\mathbf{2}, \; \mathbf{5}$ \\
\end{tabular}
\end{equation*}

\subsection{Surfaces of type $III$}
$G=\mathbb{Z}_2\times \mathbb{Z}_4, \quad \mathbf{m}=(2^2,4^2),
\quad g(C)=5$.\\  Consider the generating vector $\{g; \; h_1, \;
h_2 \}$; condition ($\ref{stab free}$) implies $g \neq (0,2)$, so up
to automorphisms of $G$ and Hurwitz moves of type $\mathbf{5}$ we
may assume $g=(1,0)$, $h_1=(0,1)$. Therefore modulo the Hurwitz move $\mathbf{1}$  we have two possibilities:
\begin{equation*}
\{(1,0); \; (0,1), \; (0,0)\} \quad \textrm{and} \quad \{(1,0); \; (0,1), \; (1,0)\},
\end{equation*}
that are equivalent via the sequence $\mathbf{1}$,
$\mathbf{3}$, $\mathbf{5}$, $\mathbf{4}$; so we may assume
$\mathcal{W}=\{(1,0); \; (0,1), \; (1,0)\}$. Now look at the
generating vector $\mathcal{V}=\{g_1, \; g_2, \; g_3, \; g_4 \}$;
here the Hurwitz moves are generated by the transposition of $g_1$
and $g_2$ and the transposition of $g_3$ and $g_4$. Condition
(\ref{stab free}) now implies $g_i \neq(1,0)$; since $\langle g_1, \ldots,
 g_4 \rangle=G$ and $g_1+g_2+g_3+g_4=0$, there are four
possibilities up to permutations:
\begin{equation*}
\begin{split}
\mathcal{V}_1 & =\{(1,2), \; (0,2), \; (0,1), \; (1,3)\}, \quad
\mathcal{V}_2=\{(1,2), \; (0,2), \; (0,3), \; (1,1)\}, \\
\mathcal{V}_3 & =\{(1,2), \; (1,2), \; (0,1), \; (0,3)\}, \quad
\mathcal{V}_4 =\{(1,2), \; (1,2), \; (1,3), \; (1,1)\}.
\end{split}
\end{equation*}
Notice that 
\begin{itemize} 
\item the automorphism of $G$ given by $(1,0) \to (1,0)$, $(0,1) \to
(0,3)$ sends $\mathcal{V}_1$ to $\mathcal{V}_2$ and $\mathcal{W}$ to
$\{(1,0); \; (0,3), \; (1,0) \}$, that is equivalent to
$\mathcal{W}$ via the Hurwitz move $\mathbf{4}$. So the pair
$(\mathcal{V}_1, \mathcal{W})$ is $\mathfrak{R}-$equivalent to
$(\mathcal{V}_2,
\mathcal{W})$; 
\item the automorphism of $G$ given by $(1,0) \to (1,0)$, $(0,1) \to
(1,3)$ sends $\mathcal{V}_3$ to $\mathcal{V}_4$ and $\mathcal{W}$ to
$\{(1,0); \; (1,3), \; (1,0) \}$, that is equivalent to
$\mathcal{W}$ via the sequence of two Hurwitz moves $\mathbf{2}$,
$\mathbf{4}$. So  $(\mathcal{V}_3, \mathcal{W})$ is
$\mathfrak{R}-$equivalent to $(\mathcal{V}_4,
\mathcal{W})$. 
\end{itemize}
On the other hand $(\mathcal{V}_1, \mathcal{W})$ and
$(\mathcal{V}_3,\mathcal{W})$ are not $\mathfrak{R}-$equivalent,
since every automorphism of $G$ leaves $(0,2)$ invariant. 
It follows that $\mathfrak{M}_{III}$ contains exactly two irreducible components.

\subsection{Surfaces of type $IV$}
$G= \mathbb{Z}_2 \times \mathbb{Z}_8, \quad \mathbf{m}=(2,8^2),
\quad g(C)=9$. \\  Consider the generating vector $\mathcal{W}=
\{g_1; \; h_1, \; h_2\}$; condition (\ref{stab free}) implies $g \neq
(0,4)$, so exactly as in the previous case we may assume, up to
automorphisms of $G$ and Hurwitz moves,
 $\mathcal{W}=\{(1,0); \; (0,1), \; (1,0)\}$. Now look at the generating
vector $\mathcal{V}=\{g_1, \; g_2, \; g_3 \}$; here the only Hurwitz
move is the transposition of $g_2$ and $g_3$.
Condition (\ref{stab free}) now implies $g_i \neq(1,0)$; since $\langle g_1,
\; g_2, \; g_3 \rangle=G$ and $g_1+g_2+g_3=0$, there are four
possibilities up to permutations:
\begin{equation*}
\begin{split}
\mathcal{V}_1 & =\{(1,4), \; (0,1), \; (1,3)\}, \quad
\mathcal{V}_2=\{(1,4), \; (1,1), \; (0,3) \}, \\
\mathcal{V}_3 & =\{(1,4), \; (1,7), \; (0,5) \}, \quad \mathcal{V}_4
=\{(1,4), \; (0,7), \; (1,5) \}.
\end{split}
\end{equation*}
Notice that
\begin{itemize}
 \item the automorphism of $G$ given by $(1,0) \to (1,0)$, $\;(0,1) \to
(1,1)$ sends $\mathcal{V}_1$ to $\mathcal{V}_2$ and $\mathcal{W}$ to
$\{(1,0); \; (1,1), \; (1,0) \}$, which is equivalent to
$\mathcal{W}$
via the Hurwitz move $\mathbf{2}$; 
\item the automorphism of $G$ given by $(1,0) \to (1,0)$, $\;(0,1) \to
(1,7)$ sends $\mathcal{V}_1$ to $\mathcal{V}_3$ and $\mathcal{W}$ to
$\{(1,0); \; (1,7), \; (1,0) \}$, which is equivalent to
$\mathcal{W}$
via the sequence of two Hurwitz moves $\mathbf{2}$, $\mathbf{4}$; 
\item the automorphism of $G$ given by $(1,0) \to (1,0)$, $\;(0,1) \to
(0,7)$ sends $\mathcal{V}_1$ to $\mathcal{V}_4$ and $\mathcal{W}$ to
$\{(1,0); \; (0,7), \; (1,0) \}$, which is equivalent to
$\mathcal{W}$
via the Hurwitz move $\mathbf{4}$. 
\end{itemize}
It follows that $(\mathcal{V}_1, \mathcal{W}), \ldots
,(\mathcal{V}_4, \mathcal{W})$ are all $\mathfrak{R}-$equivalent,
hence $\mathfrak{M}_{IV}$ is irreducible. \\ \\This completes the
proof of Theorem \ref{modulispaces}.

\section{Abelian case: the paracanonical system}  \label{paracansystem}
Now we want to study the paracanonical system of surfaces
constructed in the previous sections. We start by recalling some
definitions and results; we refer the reader to [CaCi91] for omitted
proofs and further details. Let $S$ be a minimal surface of general
type with $p_g=q=1$, let $\alpha \colon S \longrightarrow E$ be its
Albanese fibration and denote by $F_t$ the fibre of $\alpha$ over
the point $t \in E$. Moreover, define $K_S+t:=K_S+F_t-F_0$, where
$0$ is the zero element in the group structure of $E$. By
Riemann-Roch we obtain
\begin{equation*}
h^0(S,\; K_S+t)=1+h^1(S, \; K_S+t)
\end{equation*}
 for all $t \in E- \{ 0 \}$. Since $p_g=1$, by semicontinuity there
is a Zariski open set $E' \subset E$, containing $0$, such that for
any $t \in E'$ we have $h^0(S, \;K_S+t)=1$; we denote by $C_t$ the
unique curve in $|K_S+t|$. The \emph{paracanonical incidence
correspondence} is the surface $Y \subset S \times E$ which is the
schematic closure of the set $ \{(x,t) \in S \times E' \; | \; x \in
C_t \}$. Then we can define $C_t$ for any $t \in E$ as the fibre of
$Y \longrightarrow E$ over $t$, and $Y$ provides in this way a flat
family of curves on $S$, that we denote by $\{K \}$ or by $\{C_t \}$
and we call the \emph{paracanonical system} of $S$. According to
[Be88], $\{K \}$ is the irreducible component of the Hilbert scheme
of curves on $S$ algebraically equivalent to $K_S$ which dominates
$E$. Let $\mathcal{P}$ be a Poincar{\'e} sheaf on $S \times E$; then
we call $\mathcal{K}=\pi_S^*(\omega_S) \otimes \mathcal{P}$ the
paracanonical system on $S \times E$. Let $\Lambda_i:= R^i (\pi_E)_*
\mathcal{K}$. By the base change theorem, $\Lambda^0$ is an
invertible sheaf on $E$,  $\Lambda^2$ is a skyscraper sheaf of
length $1$ supported at the origin and $\Lambda^1$ is zero at the
origin, and supported on the set of points $\{ t \in E \; | \;
h^0(S, \;K_S+t) > 1 \}$; set $\lambda:= \textrm{length}(\Lambda^1)$.
\begin{definition}
The \emph{index} $\iota=\iota(K)$ of the paracanonical system is the
intersection number $Y \cdot (\{ x \} \times E)$. Roughly speaking,
$\iota$ is the number of paracanonical curves through a general
point of $S$.
\end{definition}
If $F$ is a smooth Albanese fibre of $S$, then the following
relation holds:
\begin{equation} \label{iota e g}
\iota=g(F) - \lambda.
\end{equation}
Set $V:= \alpha_*\omega_S$. Then $V$ is a vector bundle of rank
$g(F)$ over $E$, such that any locally free quotient $Q$ of $V$
verifies $\textrm{deg}(Q) \geq 0$ (this is a consequence of Fujita's
theorem, see \cite{Fu78}). Moreover we have
\begin{equation} \label{numeri V}
h^0(E, \; V)=1; \quad h^1(E, \; V)=0; \quad \textrm{deg}(V)=1.
\end{equation}
By Krull-Schmidt theorem (see \cite{At56}) there is a decomposition
of $V$ into irreducible summands:
\begin{equation} \label{decomp. irred.}
V = \bigoplus_{i=1}^k W_i
\end{equation}
which is unique up to isomorphisms. Set $d_i:= \textrm{deg}(W_i)$;
by (\ref{numeri V}) we may assume $d_1=1$, and $d_i=0$ for $2 \leq i
\leq k$. The following result shows that decomposition (\ref{decomp.
irred.}) is strongly related on the behavior of the paracanonical
system $\{ K \}$.
 \begin{proposition} \label{decomp parac}
Let $V= \bigoplus_{i=1}^k W_i$ as above. Then the following holds:
\begin{itemize}
\item[$(i)$] $k= \lambda+1$; \item[$(ii)$] \emph{rank}$(W_1)=
\iota$; \item[$(iii)$] \emph{rank}$(W_i)=1$ for $2 \leq i \leq k$.
Hence $W_i$ is a line bundle of degree $0$ for $i >1$;
\item[$(iv)$]let $L$ be a line bundle over $E$; then $h^0(S, \;
\omega_S \otimes \alpha^* L)>1$  if and only if $L=W_i^{-1}$ for
some $i>1$.
\end{itemize}
\end{proposition}
\begin{proof}
See [CaCi91].
\end{proof}
Now we can prove the main result of this section.
\begin{theorem} \label{calcolo indice}
Let $S=(C\times F)/G$ be a surface of general type with $p_g=q=1$,
isogenous to an unmixed product. If $G$ is abelian, then
$\iota(K)=1$.
\end{theorem}
\begin{proof}
We start with a lemma.
\begin{lemma} \label{aabb}
Let $\pi_C \colon C \times F \lr C, \; \pi_F \colon C \times F \lr
F$ be the two projections. Then $(\pi_C)_*\pi_F^*
\omega_F=\mathcal{O}_C^{\oplus g(F)}$.
\end{lemma}
\begin{proof}
Being $C \times F$ a product, if we fix one fibre $F_o$ of the map
$\pi_C$ then \emph{any} fibre of the bundle $(\pi_C)_* \pi_F^*
\omega_F$ can be canonically identified with the vector space
$H^0(F_o, (\pi_F^* \;  \omega_F)| _{F_o})$, which in turn is
isomorphic to $H^0(F_o, \omega_{F_o})=\mathbb{C}^{g(F)}$ by the
adjunction formula. This ends the proof. \\
\end{proof}
Now consider the commutative diagram
\begin{equation} \label{base change}
\begin{CD}
C \times F  @>p>> S\\
@VV \pi_C V     @VV \alpha V\\
C @>h>> E. \\
\end{CD}
\end{equation}
Since flatness commutes with the base change (see \cite{Ha2}), we
have
\begin{equation*}
\alpha_*p_*\omega_{C \times F} = h_*(\pi_C)_*\omega_{C \times F}.
\end{equation*}
On the other hand, by using projection formula and Lemma \ref{aabb}
we can write
\begin{equation*}
\begin{split}
(\pi_C)_*\omega_{C \times F}&=(\pi_C)_*(\pi_F^* \; \omega_F \otimes
\pi_C^*  \; \omega_C)\\ &=(\pi_C)_*\pi_F^*  \; \omega_F \otimes
\omega_C= \omega_C^{\oplus g(F)}.
\end{split}
\end{equation*}
Hence we obtain
\begin{equation} \label{relationbeta1}
\alpha_*p_*\omega_{C \times F}=(h_* \omega_C)^{\oplus g(F)}.
\end{equation}
Since $G$ is abelian, the structure theorem for abelian covers
proven in [Par91] implies that the sheaves $p_* \mathcal{O}_{C
\times F}$ and $h_* \mathcal{O}_C$ split in the following way:
\begin{equation} \label{decomp abeliana}
\begin{split}
p_* \mathcal{O}_{C \times F} & = \mathcal{O}_S \oplus
\bigoplus_{\chi \in G^* \backslash \{0\}} \mathcal{L}_{\chi}^{-1} \\
h_* \mathcal{O}_C & = \mathcal{O}_E \oplus \bigoplus_{\chi \in G^*
\backslash \{0\}} L_{\chi}^{-1},
\end{split}
\end{equation}
where $G^*$ is the group of irreducible characters of $G$ and
$\mathcal{L}_{\chi}, \; L_{\chi}$ are line bundles. More precisely,
$\mathcal{L}_{\chi}^{-1}$ and $L_{\chi}^{-1}$ are the eigensheaves
corresponding to the non-zero character $\chi \in G^*$. Moreover,
since the map $p \colon C \times F \longrightarrow S$ is {\'e}tale,
the degree of each $\mathcal{L}_{\chi}$ is zero. From (\ref{decomp
abeliana}) we obtain
\begin{equation*}
\begin{split}
p_* \omega_{C \times F} & = \omega_S \oplus  \bigoplus_{\chi \in G^*
 \backslash \{0\}}
(\omega_S \otimes \mathcal{L}_{\chi}), \\
h_* \omega_C & = \mathcal{O}_E \oplus \bigoplus_{\chi \in G^*
\backslash \{0\}} L_{\chi},
\end{split}
\end{equation*}
that is, using relation (\ref{relationbeta1}),
\begin{equation} \label{relazioneAAAAA}
\alpha_*\omega_S \oplus \bigoplus_{\chi \in G^* \backslash \{0\}}
\alpha_*(\omega_S \otimes \mathcal{L}_{\chi})= \mathcal{O}_E^{\oplus
g(F)}\bigoplus_{\chi \in G^* \backslash \{0\}}L_{\chi}^{\oplus
g(F)}.
\end{equation}
The right-hand side of (\ref{relazioneAAAAA}) is a direct sum of
line bundles; since the decomposition of a vector bundle into
irreducible summands is unique up to isomorphisms, we deduce that
$\alpha_* \omega_S$ decomposes as a direct sum of line bundles. Then
$\textrm{rank}(W_1)=1$, which implies $\iota(K)=1$ by Proposition
\ref{decomp parac} $(ii)$. This concludes the proof of Theorem
\ref{calcolo indice}.
\end{proof}
If $S$ is any minimal surface of general type with $p_g=q=1$, let us
write $\{K \}= Z+\{M \}$, where $Z$ is the fixed part and $\{M \}$
is the movable part of the paracanonical system.
\begin{corollary} \label{alb=par}
Let $S$ as in Theorem \emph{\ref{calcolo indice}}. Then  $\{M\}$
coincides with the Albanese pencil $\{F\}$.
\end{corollary}
\begin{proof}
Since $\iota=1$, through the general point of $S$ passes only one
paracanonical curve, hence $M^2=0$. By \cite{CaCi1}, Lemma 3.1 the
general member of $\{M \}$ is irreducible, hence $\{M\}$ provides an
connected, irrational pencil on $S$. By the universal property of the Albanese
morphism, it follows $\{M \}= \{F \}$. 
\end{proof}

\section{The nonabelian case} \label{newexamples}
The classification of surfaces of general type with $p_g=q=1$,
isogenous to a product of unmixed type, is still lacking when the
group $G$ is not abelian. The following theorem sheds some light on
this problem, by providing several examples.
\begin{theorem} \label{noabeliano}
Let $S=(C \times F) /G$ be a surface of general type with $p_g=q=1$,
isogenous to an unmixed product, and suppose that the group $G$ is
not abelian. Then the following cases occur.
\begin{table}[ht!]
\begin{center}
\begin{tabular}{c|c|c|c}
$G$ & $|G|$ & $g(C)$ & $g(F)$ \\
\hline 
 $S_3$ &  $6$ & $3$ & $4$  \\
 $D_4$ & $8$ & $3$ & $5$ \\
 $D_6$ & $12$ & $7$ & $3$  \\
 $A_4$ & $12$ & $4$ & $5$  \\
 $S_4$ & $24$ & $9$ & $4$ \\
 $A_5$ & $60$ & $21$ & $4$ 
\end{tabular}
\end{center}
\end{table}
\end{theorem}
The remainder of Section \ref{newexamples} deals with the
proof of Theorem \ref{noabeliano}. Let $S=(C \times F) /G$ be a
minimal surface of general type with $p_g=q=1$, isogenous to an
unmixed product, and let $f \colon F \lr
\mathbb{P}^1, \; h \colon C \lr E$ be the two quotient maps. Therefore
$f, \; h$ are induced by two admissible epimorphisms
\begin{equation*}
\vartheta  \colon \Gamma(0 \;|\;\mathbf{m}) \lr G, \quad \psi
\colon \Gamma(1 \;|\; \mathbf{n}) \lr G,
\end{equation*}
where $\mathbf{m}=(m_1, \ldots, m_r), \; \mathbf{n}=(n_1,
\ldots,n_s)$. Let $\mathcal{V}=\{g_1, \ldots,g_r \}$ and
$\mathcal{W}=\{\ell_1, \ldots, \ell_s; \; h_1,h_2 \}$ be the
generating vectors defined by $\vartheta$ and $\psi$, respectively.
By definition we have
\begin{equation*}
\begin{split}
g_1^{m_1} =\cdots =g_r^{m_r}&=g_1g_2 \cdots g_r=1,\\
 \ell_1^{n_1}=\cdots =\ell_s^{n_s}&=\ell_1\ell_2 \cdots \ell_s[h_1,h_2]=1,\\
G=\langle g_1, \ldots, g_r \rangle &=\langle \ell_1, \ldots, \ell_s,
\; h_1, \; h_2 \rangle.
\end{split}
\end{equation*}
The cyclic subgroups $\langle g_1 \rangle, \ldots ,\langle g_r
\rangle$ and their conjugates are the non-trivial stabilizers of the
action of $G$ on $F$, whereas $\langle \ell_1 \rangle, \ldots,
\langle \ell_s \rangle$ and their conjugates are the non-trivial
stabilizers of the actions of $G$ on $C$; then the diagonal action
of $G$ on $C \times F$ is free if and only if
\begin{equation} \label{stabilizzatori}
\bigg( \bigcup_{h \in G} \bigcup_{i=1}^r \langle h g_i h^{-1}
\rangle \bigg) \cap \bigg( \bigcup_{h \in G} \bigcup_{j=1}^s \langle
h \ell_j h^{-1} \rangle \bigg) =\{1 \}.
\end{equation}
Summing up, we obtain the following generalization of
 Proposition \ref{structureresult1} to the nonabelian case.
\begin{proposition} \label{structureresult2}
Let us suppose that we have the following data:
\begin{itemize}
\item a finite group $G$;
\item two admissible epimorphisms
\begin{equation*}
\begin{split}
& \vartheta \colon \Gamma(0 \; | \; \mathbf{m}) \lr G, \quad  \mathbf{m}=(m_1, \ldots, m_r)  \\
& \psi \colon \Gamma(1 \; | \; \mathbf{n}) \lr G, \quad \; \;
\mathbf{n}=(n_1, \ldots, n_s)
\end{split}
\end{equation*}
with corresponding generating vectors $\mathcal{V}=\{g_1, \ldots, g_r
\}$ and $\mathcal{W}=\{\ell_1, \ldots, \ell_s; \; h_1,h_2 \}$.
\end{itemize}
Let
\begin{equation*}
\begin{split}
& f \colon F \lr \mathbb{P}^1=F/G \\
& h \colon C \lr E= C/G
\end{split}
\end{equation*}
be the  $G-$coverings induced by $\vartheta$ and $\psi$  and let
$g(F), \; g(C)$ be the genera of $F$ and $ C$, that are related on
$|G|, \; \mathbf{m}, \; \mathbf{n}$ by $(\ref{generi1})$. Assume
moreover that
\begin{itemize}
\item $g(C) \geq 3, \;\;g(F) \geq3$;
 \item $|G|=(g(C)-1)(g(F)-1)$;
 \item condition $(\ref{stabilizzatori})$ is satisfied.
\end{itemize}
Then the diagonal action of $G$ on $C \times F$ is free and the
quotient $S=(C \times F)/G$ is a minimal surface of general type
with $p_g=q=1$. Conversely, any surface of general type with
$p_g=q=1$, isogenous to an unmixed product, arises in this way.
\end{proposition}
\begin{remark} \label{moduli-new}
We could also generalize Proposition \ref{structuremoduli} to the
nonabelian case, in order to study  the moduli spaces of surfaces
listed in Theorem \ref{noabeliano}, but we will not develop this
point here.
\end{remark}
\begin{remark} \label{achieve-non-ab}
By Proposition \ref{various-cases} we have $g(F) \leq 5$, so
$|\textrm{Aut}(G)| \leq 192$ (see \cite{Br90}, p. 91). We believe
that the classification of the unmixed, nonabelian case is not out of
reach and we hope to achieve it on a forthcoming paper.
\end{remark}

Now let us construct our examples. In the case of symmetric groups,
we will write the composition of permutations from the right to the
left; for instance, $(13)(12)=(123)$.

\subsection{$G=S_3, \;\; g(C)=3, \;\; g(F)=4$} \label{S3}
 $ $\\Take $\mathbf{m}=(2^6), \; \mathbf{n}=(3^1)$ and set
\begin{equation*}
\begin{split}
g_1&=g_2=(12), \; \; g_3=g_4=(13), \; \; g_5=g_6=(23) \\
h_1&=(12), \; \; h_2=(123), \; \; \ell_1=(132).
\end{split}
\end{equation*}
Condition (\ref{stabilizzatori}) is satisfied, hence Proposition
\ref{structureresult2} implies that this case occurs. 

\subsection{$G=D_4, \;\; g(C)=3, \;\; g(F)=5$} \label{D4}
$ $ \\$G$ is the group of order $8$ with presentation
\begin{equation*}
\langle \rho, \sigma \; | \; \rho^4=\sigma^2=1, \; \rho
\sigma=\sigma\rho^{3} \rangle.
\end{equation*}
Take $\mathbf{m}=(2^6), \; \mathbf{n}=(2^1)$ and set
\begin{equation*}
\begin{split}
g_1&=g_2=g_3=g_4=\sigma, \;\; g_5=g_6=\rho \sigma \\
h_1&=\sigma, \; \; h_2= \rho, \; \; \ell_1=\rho^2.
\end{split}
\end{equation*}
Condition (\ref{stabilizzatori}) is  satisfied, so this case occurs.
This example and the previous one were already described in \cite{Pol06}.

\subsection{$G=D_6, \;\; g(C)=7, \;\; g(F)=3$} \label{D6}
$ $ \\$G$ is the group of order $12$ with presentation
\begin{equation*}
\langle \rho, \sigma \; | \; \rho^6=\sigma^2=1, \; \rho
\sigma=\sigma\rho^{5} \rangle.
\end{equation*}
Take $\mathbf{m}=(2^3,6^1), \; \mathbf{n}=(2^2)$ and set
\begin{equation*}
\begin{split}
g_1&=\rho^3, \; \; g_2=\rho \sigma, \;\; g_3=\rho^5 \sigma, \; \;
g_4=\rho \\
h_1&=h_2= \rho, \;\; \ell_1=\ell_2=\sigma.
\end{split}
\end{equation*}
Condition (\ref{stabilizzatori}) is satisfied, so
 this case occurs.

\subsection{$G=A_4, \;\; g(C)=4, \;\; g(F)=5$} \label{A4}
$ $\\Take $\mathbf{m}=(3^4), \; \mathbf{n}=(2^1)$ and set
\begin{equation*}
\begin{split}
g_1&=(234), \; \; g_2=(123), \; \; g_3=(124), \; \; g_4=(134) \\
h_1&=(123), \;\; h_2=(124), \;\; \ell_1=(12)(34).
\end{split}
\end{equation*}
Condition (\ref{stabilizzatori}) is  satisfied, so this case occurs.

\subsection{$G=S_4, \; \; g(C)=9, \; \; g(F)=4$} \label{S4}
$ $ \\Take $\mathbf{m}=(2^3,4^1), \; \mathbf{n}=(3^1)$ and set
\begin{equation*}
\begin{split}
g_1&=(23), \; \; g_2=(24), \;\; g_3=(12), \;\; g_4=(1234)\\
h_1&=(12), \;\; h_2=(1234), \;\; \ell_1=(132).
\end{split}
\end{equation*}
Condition (\ref{stabilizzatori}) is satisfied, so this case occurs.

\subsection{$G=A_5, \; \; g(C)=21, \;\; g(F)=4$} \label{A5}
$ $\\Take $\mathbf{m}=(2,5^2), \; \mathbf{n}=(3^1)$ and set
\begin{equation*}
\begin{split}
g_1&=(24)(35), \; \; g_2=(13452), \;\; g_3=(12345) \\
h_1&=(345), \;\; h_2=(15432), \;\; \ell_1=(235).
\end{split}
\end{equation*}
One checks by direct computation that $g_1g_2g_3=\ell_1[h_1,
h_2]=1$. Since $g_3g_1g_3=(152)$, it follows that the subgroup
generated by $g_1,g_2,g_3$ has order at least $2 \cdot 3 \cdot
5=30$. On the other hand $G$ is simple, so it cannot contain a
subgroup of order $30$; therefore $\langle g_1,g_2,g_3 \rangle =G$.
Analogously,  $\ell_1 h_1 \ell_1=(24)(35)$ which implies that the
subgroup $\langle \ell_1, h_1, h_2 \rangle$
 has order at least $30$ and so must be equal to $G$ too.
Condition (\ref{stabilizzatori}) is verified, hence this case occurs. \\
\\ This completes the proof of Theorem \ref{noabeliano}.

\end{document}